\documentclass[reqno]{amsart}
\usepackage{setspace}
\usepackage{graphicx}
\usepackage{amsmath}
\usepackage{amsfonts}
\usepackage{amssymb}

\begin{document}

\centerline{\huge \bf Characterization of the Two-Dimensional}

\medskip
\centerline{\huge \bf Five-Fold Lattice Tiles}

\bigskip
\centerline{\large Chuanming Zong}

\vspace{0.6cm}
\centerline{\begin{minipage}{12.5cm}
{\bf Abstract.} In 1885, Fedorov discovered that a convex domain can form a lattice tiling of the Euclidean plane if and only if it is a parallelogram or a centrally symmetric hexagon. It is known that there is no other convex domain which can form a two-, three- or four-fold lattice tiling in the Euclidean plane, but there are centrally symmetric convex octagons and decagons which can form five-fold lattice tilings. This paper characterizes all the convex domains which can form five-fold lattice tilings of the Euclidean plane. They are parallelograms, centrally symmetric hexagons, centrally symmetric octagons (under suitable affine linear transformations) with vertices ${\bf v}_1=(-\alpha , -{3\over 2})$, ${\bf v}_2=(1-\alpha , -{3\over 2})$, ${\bf v}_3=(1+\alpha , -{1\over 2})$, ${\bf v}_4=(1-\alpha , {1\over 2})$, ${\bf v}_5=-{\bf v}_1$, ${\bf v}_6=-{\bf v}_2$, ${\bf v}_7=-{\bf v}_3$ and ${\bf v}_8=-{\bf v}_4$, where $0<\alpha <{1\over 4}$, or with vertices ${\bf v}_1=(\beta , -2)$, ${\bf v}_2=(1+\beta , -2)$, ${\bf v}_3=(1-\beta , 0)$, ${\bf v}_4=(\beta , 1)$, ${\bf v}_5=-{\bf v}_1$, ${\bf v}_6=-{\bf v}_2$, ${\bf v}_7=-{\bf v}_3$, ${\bf v}_8=-{\bf v}_4$, where ${1\over 4}<\beta <{1\over 3}$, or centrally symmetric decagons (under suitable affine linear transformations) with ${\bf u}_1=(0, 1)$, ${\bf u}_2=(1, 1)$, ${\bf u}_3=({3\over 2}, {1\over 2})$, ${\bf u}_4=({3\over 2}, 0)$, ${\bf u}_5=( 1,-{1\over 2})$, ${\bf u}_6=-{\bf u}_1$, ${\bf u}_7=-{\bf u}_2$, ${\bf u}_8=-{\bf u}_3$, ${\bf u}_9=-{\bf u}_4$ and ${\bf u}_{10}=-{\bf u}_5$ as the middle points of their edges.
\end{minipage}}


\vspace{0.8cm}
\noindent
{\Large\bf 1. Introduction}

\bigskip\noindent
Planar tilings is an ancient subject in our civilization. It has been considered in the arts by craftsmen since antiquity. Up to now, it is still an active research field in mathematics and some basic problems remain unsolved. In 1885, Fedorov \cite{fedo} discovered that there are only two types of two-dimensional lattice tiles: {\it parallelograms and centrally symmetric hexagons}. In 1917, for the purpose to verify the second part of Hilbert's 18th problem in $\mathbb{E}^2$, Bieberbach suggested Reinhardt (see \cite{rein}) to determine all the two-dimensional congruent tiles. However, to complete the list turns out to be challenging and dramatic. Over the years, the list has been successively extended by Reinhardt, Kershner, James, Rice, Stein, Mann, McLoud-Mann and Von Derau (see \cite{mann,zong14}), its completeness has been mistakenly announced several times! In 2017, M. Rao \cite{rao} announced a completeness proof based on computer checks.

The three-dimensional case was also studied in the ancient time. More than 2,300 years ago, Aristotle claimed that both identical regular tetrahedra and identical cubes can fill the whole space without gap. The cube case is obvious! However, the tetrahedron case is wrong and such a tiling is impossible (see \cite{lazo}).

Let $K$ be a convex body with (relative) interior ${\rm int}(K)$, (relative) boundary $\partial (K)$ and volume ${\rm vol}(K)$, and let $X$ be a discrete set, both in $\mathbb{E}^n$. We call $K+X$ a {\it translative tiling} of $\mathbb{E}^n$ and call $K$ a {\it translative tile} if $K+X=\mathbb{E}^n$ and the translates ${\rm int}(K)+{\bf x}_i$ are pairwise disjoint. In other words, if $K+X$ is both a packing and a covering in $\mathbb{E}^n$. In particular, we call $K+\Lambda$ a {\it lattice tiling} of $\mathbb{E}^n$ and call $K$ a {\it lattice tile} if $\Lambda $ is an $n$-dimensional lattice. Apparently, a translative tile must be a convex polytope. Usually, a lattice tile is called a {\it parallelohedron}.

In 1885, Fedorov \cite{fedo} also characterized the three-dimensional lattice tiles: {\it A three-dimensional lattice tile must be a parallelotope, an hexagonal prism, a rhombic dodecahedron, an elongated dodecahedron, or a truncated octahedron.} The situations in higher dimensions turn out to be very complicated. Through the works of Delone \cite{delo}, $\check{S}$togrin \cite{stog} and Engel \cite{enge}, we know that there are exact $52$ combinatorially different types of parallelohedra in $\mathbb{E}^4$. A computer classification for the five-dimensional parallelohedra was announced by Dutour Sikiri$\acute{\rm c}$, Garber, Sch$\ddot{\rm u}$rmann and Waldmann \cite{dgsw} only in 2015.

Let $\Lambda $ be an $n$-dimensional lattice. The {\it Dirichlet-Voronoi cell} of $\Lambda $ is defined by
$$C=\left\{ {\bf x}: {\bf x}\in \mathbb{E}^n,\ \| {\bf x}, {\bf o}\|\le \| {\bf x}, \Lambda \|\right\},$$
where $\| X, Y\|$ denotes the Euclidean distance between $X$ and $Y$. Clearly, $C+\Lambda $ is a lattice tiling and the Dirichlet-Voronoi cell $C$ is a parallelohedron. In 1908,
Voronoi \cite{voro} made a conjecture that {\it every parallelohedron is a linear transformation image of the Dirichlet-Voronoi cell of a suitable lattice.} In $\mathbb{E}^2$, $\mathbb{E}^3$ and $\mathbb{E}^4$, this conjecture was confirmed by Delone \cite{delo} in 1929. In higher dimensions, it is still open.

To characterize the translative tiles is another fascinating problem. At the first glance, translative tilings should be more complicated than lattice tilings. However, the dramatic story had a happy end! It was shown by Minkowski \cite{mink} in 1897 that {\it every translative tile must be centrally symmetric}. In 1954, Venkov \cite{venk} proved that {\it every translative tile must be a lattice tile $($parallelohedron$)$} (see \cite{alek} for generalizations). Later, a new proof for this beautiful result was independently discovered by McMullen \cite{mcmu}.

Let $X$ be a discrete multiset in $\mathbb{E}^n$ and let $k$ be a positive integer. We call $K+X$ a {\it $k$-fold translative tiling} of $\mathbb{E}^n$ and call $K$ a {\it $k$-fold translative tile} if every point ${\bf x}\in \mathbb{E}^n$ belongs to at least $k$ translates of $K$ in $K+X$ and every point ${\bf x}\in \mathbb{E}^n$ belongs to at most $k$ translates of ${\rm int}(K)$ in ${\rm int}(K)+X$. In other words, $K+X$ is both a $k$-fold packing and a $k$-fold covering in $\mathbb{E}^n$. In particular, we call $K+\Lambda$ a {$k$-fold lattice tiling} of $\mathbb{E}^n$ and call $K$ a {\it $k$-fold lattice tile} if $\Lambda $ is an $n$-dimensional lattice. Apparently, a $k$-fold translative tile must be a convex polytope. In fact, similar to Minkowski's characterization, it was shown by Gravin, Robins and Shiryaev \cite{grs} that {\it a $k$-fold translative tile must be a centrally symmetric polytope with centrally symmetric facets.}

Multiple tilings was first investigated by Furtw\"angler \cite{furt} in 1936 as a generalization of Minkowski's conjecture on cube tilings. Let $C$ denote the $n$-dimensional unit cube. Furtw\"angler made a conjecture that {\it every $k$-fold lattice tiling $C+\Lambda$ has twin cubes. In other words, every multiple lattice tiling $C+\Lambda$ has two cubes sharing a whole facet.} In the same paper, he proved the two- and three-dimensional cases. Unfortunately, when $n\ge 4$, this beautiful conjecture was disproved by Haj\'os \cite{hajo} in 1941. In 1979, Robinson \cite{robi} determined all the integer pairs $\{ n,k\}$ for which Furtw\"angler's conjecture is false. We refer to Zong \cite{zong05,zong06} for an introduction account and a detailed account on this fascinating problem, respectively, to pages 82-84 of Gruber and Lekkerkerker \cite{grub} for some generalizations.

Let $P$ be an $n$-dimensional centrally symmetric convex polytope, let $\tau (P)$ denote the smallest integer $k$ such that $P$ is a $k$-fold translative tile, and let $\tau^* (P)$ denote the smallest integer $k$ such that $P$ is a $k$-fold lattice tile. For convenience, we define $\tau (P)=\infty $ if $P$ can not form translative tiling of any multiplicity. Clearly, for every convex polytope we have
$$\tau (P)\le \tau^*(P).$$
If $\sigma $ is a non-singular affine linear transformation from $\mathbb{E}^n$ to $\mathbb{E}^n$, it can be easily verified that $P+X$ is a $k$-fold tiling of $\mathbb{E}^n$ if and only if $\sigma (P)+\sigma (X)$ is a $k$-fold tiling of $\mathbb{E}^n$. Thus, both $\tau (\sigma (P))=\tau (P)$ and $\tau^*(\sigma (P))=\tau^*(P)$ hold for all convex polytopes $P$ and all non-singular affine linear transformations $\sigma$.

In 1994, Bolle \cite{boll} proved that {\it every centrally symmetric lattice polygon is a multiple lattice tile}. However, little is known about the multiplicity. Let $\Lambda $ denote the two-dimensional integer lattice, and let $D_8$ denote the octagon with vertices $(1,0)$, $(2,0)$, $(3,1)$, $(3,2)$, $(2,3)$, $(1,3)$, $(0,2)$ and $(0,1)$. As a particular example of Bolle's theorem, it was discovered by Gravin, Robins and Shiryaev \cite{grs} that {\it $D_8+\Lambda$ is a seven-fold lattice tiling of $\mathbb{E}^2$.} Consequently, we have
$$\tau^*(D_8)\le 7.$$

In 2000, Kolountzakis \cite{kolo} proved that, if $D$ is a two-dimensional convex domain which is not a parallelogram and $D+X$ is a multiple tiling in $\mathbb{E}^2$, then $X$ must be a finite union of translated two-dimensional lattices. In 2013, a similar result in $\mathbb{E}^3$ was discovered by Gravin, Kolountzakis, Robins and Shiryaev \cite{gkrs}.

\medskip
Recently, Yang and Zong \cite{yz1} proved the following results: {\it Besides parallelograms and centrally symmetric hexagons, there is no other convex domain which can form a two-, three- or four-fold lattice tiling in the Euclidean plane. However, there are convex octagons and decagons which can form five-fold lattice tilings. Consequently, whenever $n\ge 3$, there are non-parallelohedral polytopes which can form five-fold lattice tilings in the $n$-dimensional Euclidean space.}

\smallskip
This paper characterizes all the two-dimensional five-fold lattice tiles by proving the following results.

\medskip\noindent
{\bf Theorem 1.} {\it A convex domain can form a five-fold lattice tiling of the Euclidean plane if and only if it is a parallelogram, a centrally symmetric hexagon, a centrally symmetric octagon (under a suitable affine linear transformation) with vertices ${\bf v}_1=(-\alpha , -{3\over 2})$, ${\bf v}_2=(1-\alpha , -{3\over 2})$, ${\bf v}_3=(1+\alpha , -{1\over 2})$, ${\bf v}_4=(1-\alpha , {1\over 2})$, ${\bf v}_5=-{\bf v}_1$, ${\bf v}_6=-{\bf v}_2$, ${\bf v}_7=-{\bf v}_3$ and ${\bf v}_8=-{\bf v}_4$, where $0<\alpha <{1\over 4},$ or with vertices
${\bf v}_1=(\beta , -2)$, ${\bf v}_2=(1+\beta , -2)$, ${\bf v}_3=(1-\beta , 0)$, ${\bf v}_4=(\beta , 1)$, ${\bf v}_5=-{\bf v}_1$, ${\bf v}_6=-{\bf v}_2$, ${\bf v}_7=-{\bf v}_3$, ${\bf v}_8=-{\bf v}_4$, where ${1\over 4}<\beta <{1\over 3}$, or a centrally symmetric decagon (under a suitable affine linear transformation) with ${\bf u}_1=(0, 1)$, ${\bf u}_2=(1, 1)$, ${\bf u}_3=({3\over 2}, {1\over 2})$, ${\bf u}_4=({3\over 2}, 0)$, ${\bf u}_5=( 1,-{1\over 2})$, ${\bf u}_6=-{\bf u}_1$, ${\bf u}_7=-{\bf u}_2$, ${\bf u}_8=-{\bf u}_3$, ${\bf u}_9=-{\bf u}_4$ and ${\bf u}_{10}=-{\bf u}_5$ as the middle points of its edges.}

\medskip\noindent
{\bf Theorem 2.} {\it Let $W$ denote the quadrilateral with vertices ${\bf w}_1=(-{1\over 2}, 1)$, ${\bf w}_2=(-{1\over 2}, {3\over 4})$, ${\bf w}_3=(-{2\over 3}, {2\over 3})$ and ${\bf w}_4=(-{3\over 4}, {3\over 4})$. A centrally symmetric convex decagon with ${\bf u}_1=(0,1)$, ${\bf u}_2=(1,1)$, ${\bf u}_3=({3\over 2}, {1\over 2})$, ${\bf u}_4=({3\over 2}, 0)$, ${\bf u}_5=(1,-{1\over 2})$, ${\bf u}_6=-{\bf u}_1$, ${\bf u}_7=-{\bf u}_2$, ${\bf u}_8=-{\bf u}_3$, ${\bf u}_9=-{\bf u}_4$ and ${\bf u}_{10}=-{\bf u}_5$ as the middle points of its edges if and only if one of its vertices is an interior point of $W$.}

\vspace{0.6cm}
\noindent
{\Large\bf 2. Basic Results}

\bigskip\noindent
Let $P_{2m}$ denote a centrally symmetric convex $2m$-gon centered at the origin, let ${\bf v}_1$, ${\bf v}_2$, $\ldots$, ${\bf v}_{2m}$ be the $2m$ vertices of $P_{2m}$ enumerated in the clock order, and let $G_1$, $G_2$, $\ldots $, $G_{2m}$ be the $2m$ edges of $P_{2m}$, where $G_i$ has two ends ${\bf v}_i$ and ${\bf v}_{i+1}$. For convenience, we write $V=\{{\bf v}_1, {\bf v}_2, \ldots, {\bf v}_{2m}\}$ and $\Gamma=\{G_1, G_2, \ldots, G_{2m}\}$.

Assume that $P_{2m}+X$ is a $\tau (P_{2m})$-fold translative tiling of $\mathbb{E}^2$, where $X=\{{\bf x}_1, {\bf x}_2, {\bf x}_3, \ldots \}$ is a discrete multiset with ${\bf x}_1={\bf o}$. Now, let us observe the local structure of $P_{2m}+X$ at the vertices ${\bf v}\in V+X$.

Let $X^{\bf v}$ denote the subset of $X$ consisting of all points ${\bf x}_i$ such that
$${\bf v}\in \partial (P_{2m})+{\bf x}_i.$$
Since $P_{2m}+X$ is a multiple tiling, the set $X^{\bf v}$ can be divided into disjoint subsets $X^{\bf v}_1$, $X^{\bf v}_2$, $\ldots ,$ $X^{\bf v}_t$ such that the translates in $P_{2m}+X^{\bf v}_j$ can be re-enumerated as $P_{2m}+{\bf x}^j_1$, $P_{2m}+{\bf x}^j_2$, $\ldots $, $P_{2m}+{\bf x}^j_{s_j}$ satisfying the following conditions:

\medskip
\noindent
{\bf 1.} {\it ${\bf v}\in \partial (P_{2m})+{\bf x}^j_i$ holds for all $i=1, 2, \ldots, s_j.$}

\smallskip\noindent
{\bf 2.} {\it Let $\angle^j_i$ denote the inner angle of $P_{2m}+{\bf x}^j_i$ at ${\bf v}$ with two half-line edges $L^j_{i,1}$ and $L^j_{i,2}$, where $L^j_{i,1}$, ${\bf x}^j_i-{\bf v}$ and $L^j_{i,2}$ are in clock order. Then, the inner angles join properly as
$$L^j_{i,2}=L^j_{i+1,1}$$
holds for all $i=1,$ $2,$ $\ldots ,$ $s_j$, where $L^j_{s_j+1,1}=L^j_{1,1}$.}

\medskip
For convenience, we call such a sequence $P_{2m}+{\bf x}^j_1$, $P_{2m}+{\bf x}^j_2$, $\ldots $, $P_{2m}+{\bf x}^j_{s_j}$ an {\it adjacent wheel} at ${\bf v}$. It is easy to see that
$$\sum_{i=1}^{s_j}\angle^j_i =2w_j\cdot \pi$$
hold for positive integers $w_j$. Then we define
$$\phi ({\bf v})=\sum_{j=1}^tw_j= {1\over {2\pi }}\sum_{j=1}^t\sum_{i=1}^{s_j}\angle^j_i$$
and
$$\varphi ({\bf v})=\sharp \left\{ {\bf x}_i:\ {\bf x}_i\in X,\ {\bf v}\in {\rm int}(P_{2m})+{\bf x}_i\right\}.$$

\medskip
Clearly, if $P_{2m}+X$ is a $\tau (P_{2m})$-fold translative tiling of $\mathbb{E}^2$, then
$$\tau (P_{2m})= \varphi ({\bf v})+\phi ({\bf v})\eqno (1)$$
holds for all ${\bf v}\in V+X$.

\medskip
First, let us introduce some basic results which will be useful in this paper.

\medskip\noindent
{\bf Lemma 1 (Bolle \cite{boll}).} {\it A convex polygon is a $k$-fold lattice tile for a lattice $\Lambda$ and some positive integer $k$ if and only if the following conditions are satisfied:

\noindent
{\bf 1.} It is centrally symmetric.

\noindent
{\bf 2.} When it is centered at the origin, in the relative interior of each edge $G$ there is a point of ${1\over 2}\Lambda $.

\noindent
{\bf 3.} If the midpoint of $G$ is not in ${1\over 2}\Lambda $ then $G$ is a lattice vector of $\Lambda $.}

\medskip\noindent
{\bf Lemma 2.} {\it If $m$ is even and $P_{2m}+\Lambda $ is a multiple lattice tiling, then $P_{2m}$ has an edge $G$ which is a lattice vector of $\Lambda $.}

\medskip\noindent
{\bf Proof.} We assume that $\Lambda =\mathbb{Z}^2$. Let ${\bf v}_1$, ${\bf v}_2$, $\ldots $, ${\bf v}_{2m}$ be the $2m$ vertices of $P_{2m}$ arranged in an anti-clock order. Let $G_i$ denote the edge of $P_{2m}$ with vertices ${\bf v}_i$ and ${\bf v}_{i+1}$, where ${\bf v}_{2m+1}={\bf v}_1$.

If the midpoint of one of the $2m$ edges, say $G_1$, is not in ${1\over 2}\Lambda $, then it follows from Lemma 1 that $G_1$ is a lattice vector of $\Lambda $.

Let ${\bf u}_i$ denote the midpoint of $G_i$. If ${\bf u}_i\in {1\over 2}\Lambda $ hold for all $i=1, 2, \ldots, 2m$, then we have
$$\left\{\begin{array}{ll}
{\bf v}_2-{\bf u}_1={\bf u}_1-{\bf v}_1, &\\
{\bf v}_3-{\bf u}_2={\bf u}_2-{\bf v}_2, &\\
\hspace{1.2cm}\ldots &\\
{\bf v}_{m+1}-{\bf u}_m={\bf u}_m-{\bf v}_m,&
\end{array}\right.$$
which implies that
$${\bf v}_{m+1}=(-1)^m{\bf v}_1+2\sum_{i=1}^m(-1)^i{\bf u}_i.\eqno (2)$$
Since $m$ is even and ${\bf v}_{m+1}=-{\bf v}_1$, it can be deduced by (2) that
$${\bf v}_1=\sum_{i=1}^m(-1)^{i+1}{\bf u}_i\in \mbox{$1\over 2$}\Lambda.$$
If fact, in this case all the vertices belong to ${1\over 2}\Lambda $. Then, we get
$${\bf v}_2-{\bf v}_1=2\left({\bf u}_1-{\bf v}_1\right) \in \Lambda .$$
The lemma is proved. \hfill{$\Box$}

\medskip\noindent
{\bf Lemma 3.} {\it Let ${\bf u}_i$ be the middle point of $G_i$. If $m$ is an odd positive integer, $P_{2m}+\Lambda $ is a $k$-fold lattice tiling of $\mathbb{E}^2$, and all ${\bf u}_i$ belong to ${1\over 2}\Lambda $, then we have
$$\sum_{i=1}^m(-1)^i{\bf u}_i ={\bf o},$$
where ${\bf o}=(0,0)$ is the origin of $\mathbb{E}^2$.}

\medskip\noindent
{\bf Proof.} Since ${\bf u}_i$ is the middle point of $G_i$, we have
$$\left\{\begin{array}{cl}
{\bf v}_2&\hspace{-0.3cm}=2{\bf u}_1-{\bf v}_1,\\
{\bf v}_3&\hspace{-0.3cm}=2{\bf u}_2-{\bf v}_2,\\
&\ldots \\
{\bf v}_{m+1}&\hspace{-0.3cm}=2{\bf u}_m-{\bf v}_m,
\end{array}\right.$$
which implies
$$-{\bf v}_1={\bf v}_{m+1}=-{\bf v}_1 - 2 \sum_{i=1}^m(-1)^i{\bf u}_i$$
and finally
$$\sum_{i=1}^m(-1)^i{\bf u}_i ={\bf o}.$$
The Lemma is proved. \hfill{$\Box$}

\medskip\noindent
{\bf Lemma 4 (Yang and Zong \cite{yz2}).} {\it Assume that $P_{2m}$ is a centrally symmetric convex $2m$-gon centered at the origin and $P_{2m}+X$ is a $\tau (P_{2m})$-fold translative tiling of the plane, where $m\ge 4$. If ${\bf v}\in V+X$ is a vertex and $G\in \Gamma +X$ is an edge with ${\bf v}$ as one of its two ends, then there are at least $\lceil (m-3)/2\rceil $ different translates $P_{2m}+{\bf x}_i$ satisfying both
$${\bf v}\in \partial (P_{2m})+{\bf x}_i$$
and}
$$G\setminus \{ {\bf v}\}\subset {\rm int}(P_{2m})+{\bf x}_i.$$

\medskip\noindent
{\bf Lemma 5 (Yang and Zong \cite{yz2}).} {\it Let $P_{2m}$ be a centrally symmetric convex $2m$-gon, then}
$$\tau^*(P_{2m})\ge \tau (P_{2m})\ge \left\{\begin{array}{ll}
m-1,&\mbox{if $m$ is even,}\\
m-2,&\mbox{if $m$ is odd.}
\end{array}
\right.$$

\vspace{0.6cm}
\noindent
{\Large\bf 3. Technical Lemmas}

\bigskip\noindent
{\bf Lemma 6.} {\it Let $P_{14}$ be a centrally symmetric convex tetradecagon, then}
$$\tau^* (P_{14})\ge \tau (P_{14})\ge 6.$$

\medskip
\noindent
{\bf Proof.} We take ${\bf v}\in V+X$ and assume that $P_{14}+{\bf x}_1$, $P_{14}+{\bf x}_2$, $\ldots $, $P_{14}+{\bf x}_s$ is an adjacent wheel at ${\bf v}$ with corresponding angles $\angle_1$, $\angle_2$, $\ldots $, $\angle_s$, where $\angle_1<\pi$. Without loss of generality, we assume further that $P_{14}+{\bf x}_1$, $P_{14}+{\bf x}_2$, $\ldots $, $P_{14}+{\bf x}_n$ is a part of the wheel such that $\angle_1$, $\angle_2$, $\ldots $, $\angle_n$ has no opposite pair, $\angle_1<\pi $ and
$$\sum_{i=1}^n\angle_i=\mu\cdot \pi, \eqno(3)$$
where $\mu $ is a positive integer.

Clearly, $\angle_i=\pi $ if and only if ${\bf v}$ is a relative interior point of an edge of $P_{14}+{\bf x}_i$ and therefore
$$\sum_{i=1}^n\angle_i<n\cdot \pi.\eqno(4)$$
On the other hand, if $\ell $ of the $n$ angles are $\pi $ and $n-\ell <m$, then we have
$$\sum_{i=1}^n\angle_i>\ell \cdot \pi +(m-1)\cdot \pi -(m-n+\ell )\cdot \pi =(n-1)\cdot \pi, \eqno(5)$$
which together with (4) contradicts (3). Therefore, to avoid the contradiction, we must have
$$n-\ell =m$$
and each pair of the opposite angles of $P_{14}$ has a representative in the angle sequence $\angle_1$, $\angle_2$, $\ldots $, $\angle_n$. Therefore, we have
$$\phi ({\bf v})\ge {1\over {2\pi }}\sum_{i=1}^s \angle_i \ge {1\over {2\pi }}\sum_{i=1}^n \angle_i\ge {1\over {2\pi }}\cdot {{(14-2)\cdot \pi }\over 2}= 3.\eqno(6)$$
On the other hand, by Lemma 4, we have
$$\varphi ({\bf v})\ge \left\lceil {{7-3}\over 2} \right\rceil = 2.\eqno(7)$$

Now, we consider two cases.

\smallskip\noindent
{\bf Case 1.} {\it $\phi ({\bf v})\ge 4$ holds for a vertex ${\bf v}\in V+X$.} Then, by (1) and (7) we get
$$\tau (P_{14})= \varphi ({\bf v}) + \phi ({\bf v}) \ge 6.$$

\smallskip\noindent
{\bf Case 2.} {\it $\phi ({\bf v})= 3$ holds for every vertex ${\bf v}\in V+X$.} First, let's observe a simple fact. If $\phi ({\bf v})=3$ holds at ${\bf v}\in V+X$ and $P_{14}+{\bf x}_1$, $P_{14}+{\bf x}_2$, $\ldots $, $P_{14}+{\bf x}_s$ is an adjacent wheel at ${\bf v}$, then it follows from (6) that $s$ must be seven and ${\bf v}$ is a common vertex of all these translates, as shown by Figure 1. Then, by Lemma 4, every vertex ${\bf v}^*_i$ connecting with ${\bf v}$ by an edge is an interior point of two of the seven translates in the wheel.

\begin{figure}[!ht]
\centering
\includegraphics[scale=0.5]{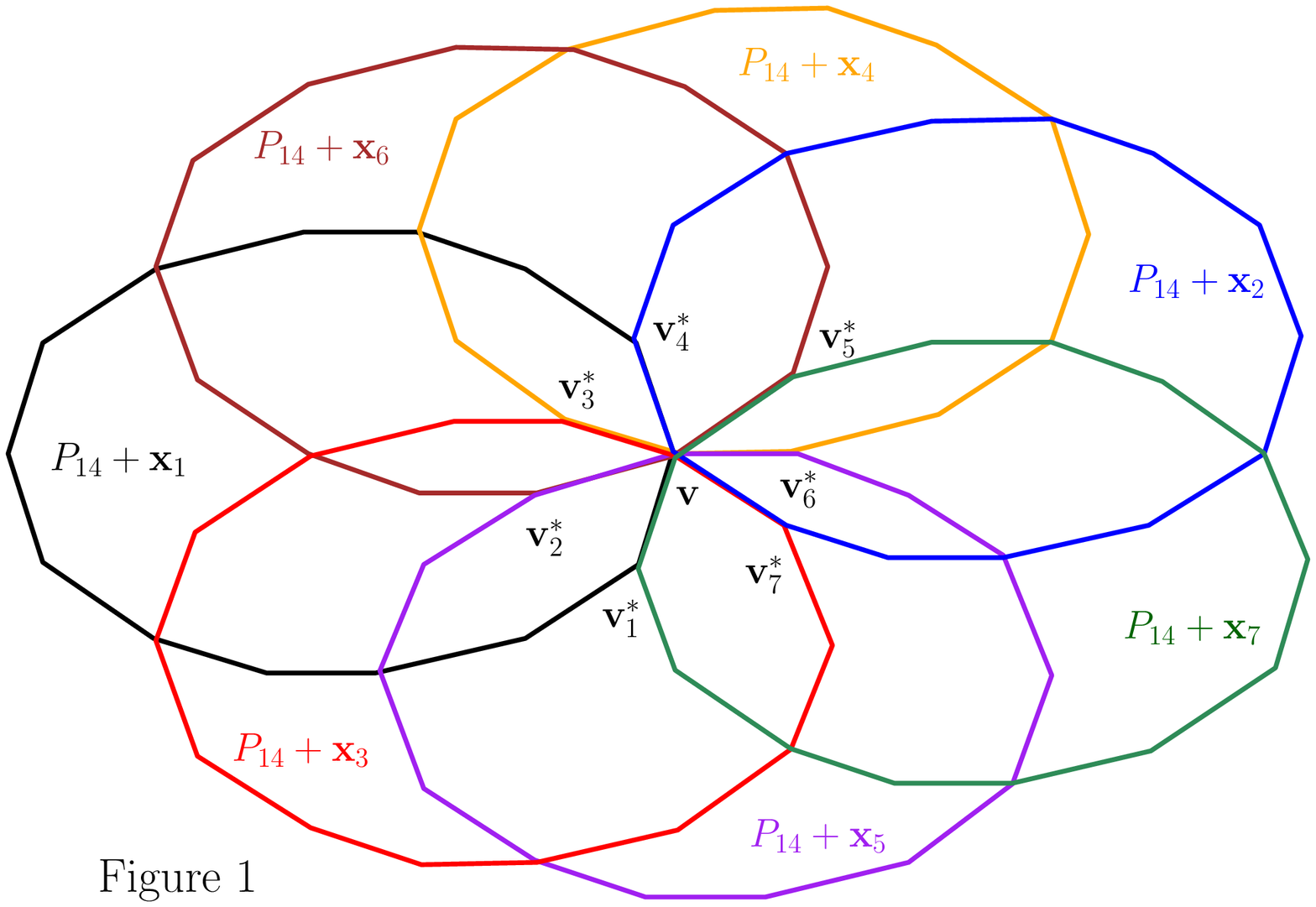}
\end{figure}

Then, by the assumption we have
$$\phi ({\bf v}^*_1)=\phi ({\bf v}^*_2)= \phi ({\bf v}^*_3)=\phi ({\bf v}^*_4)=\phi ({\bf v}^*_5)=\phi ({\bf v}^*_6)=\phi ({\bf v}^*_7)=3.$$
Therefore, for each vertex ${\bf v}^*_i$ there are two different points ${\bf y}_{i,1}$, ${\bf y}_{i,2}\in X$ such that
$${\bf v}^*_i\in \partial (P_{14})+{\bf y}_{i,j},\qquad j=1,\ 2$$
and
$${\bf v}\in {\rm int}(P_{14})+{\bf y}_{i,j},\qquad j=1,\ 2.$$

If ${\bf y}_{i,j}\not\in \{{\bf y}_{1,1}, {\bf y}_{1,2}\}$ holds for one of these points, then we get
$$\varphi ({\bf v})\ge 3$$
and therefore
$$\tau (P_{14})\ge \varphi ({\bf v}) + \phi ({\bf v}) \ge 6.$$

If ${\bf y}_{i,j}\in \{{\bf y}_{1,1}, {\bf y}_{1,2}\}$ holds for all of these points, then we must have
$$\{{\bf v}^*_1, {\bf v}^*_2, {\bf v}^*_3, {\bf v}^*_4, {\bf v}^*_5, {\bf v}^*_6, {\bf v}^*_7\} \subset \partial (P_{14})+{\bf y}_{1,1}.$$
It is known that $(D+{\bf x})\cap (D+{\bf y})$ is centrally symmetric for all ${\bf x}$ and ${\bf y}$ whenever $D$ is centrally symmetric. Then, by Figure 1 it is easy to see that $(P_{14}+{\bf y}_{1,1})\cap (P_{14}+{\bf x}_1)$ is a parallelogram with vertices ${\bf v}^*_1$, ${\bf v}$, ${\bf v}^*_4$ and ${\bf v}^*_1+({\bf v}^*_4-{\bf v})$, and $(P_{14}+{\bf y}_{1,1})\cap (P_{14}+{\bf x}_7)$ is a parallelogram with vertices ${\bf v}^*_1$, ${\bf v}$, ${\bf v}^*_5$ and ${\bf v}^*_1+({\bf v}^*_5-{\bf v})$. Consequently, by symmetry one can deduce that $P_{14}+{\bf y}_{1,1}$ is an hexagon with vertices ${\bf v}^*_1$, ${\bf v}^*_1+({\bf v}^*_4-{\bf v})$, ${\bf v}^*_4$, ${\bf v}+({\bf v}-{\bf v}^*_1)$, ${\bf v}^*_5$ and ${\bf v}^*_1+({\bf v}^*_5-{\bf v})$, which contradicts the assumption that $P_{14}$ is a tetradecagon.

As a conclusion, for every centrally symmetric convex tetradecagon we have
$$\tau (P_{14})\ge 6.$$
The lemma is proved.\hfill{$\Box$}

\medskip
\noindent
{\bf Lemma 7.} {\it Let $P_{12}$ be a centrally symmetric convex dodecagon, then}
$$\tau^* (P_{12})\ge 6.$$

\medskip
\noindent
{\bf Proof.} Since $\tau^*(P_{2m})$ is invariant under linear transformations on $P_{2m}$, we assume that $\Lambda =\mathbb{Z}^2$ and $P_{12}+\Lambda $ is a $\tau^*(P_{12})$-fold lattice tiling. Let ${\bf u}_i$ denote the middle point of $G_i$ and write ${\bf v}_i=(x_i,y_i)$ and ${\bf u}_i=(x'_i, y'_i)$. By Lemma 2 and a uni-modular transformation, we may assume that ${\bf v}_2-{\bf v}_1=(k, 0)$ and $y'_1>0$, where $k$ is a positive integer. By reduction, we many assume further that ${\bf v}_2-{\bf v}_1=(1,0)$. For convenience, let $P$ denote the parallelogram with vertices ${\bf v}_1$, ${\bf v}_2$, ${\bf v}_7=-{\bf v}_1$
and ${\bf v}_8=-{\bf v}_2$.

By Lemma 1 it follows that all $y_2-y_3$, $y_3-y_4$, $y_4-y_5$, $y_5-y_6$ and $y_6-y_7$ are positive integers. Thus, we have
$$y_1=y'_1=y_2\ge {5\over 2},$$
$${\rm vol} (P_{12})>{\rm vol} (P)\ge 5$$
and therefore, since $\tau^*(P_{12})$ is an integer,
$$\tau^*(P_{12})={\rm vol} (P_{12})\ge 6.$$
The lemma is proved. \hfill{$\Box$}

\medskip\noindent
{\bf Lemma 8.} {\it For every centrally symmetric convex decagon $P_{10}$ we have
$$\tau^*(P_{10})\ge 5,$$
where the equality holds if and only if, under a suitable affine linear transformation, it takes ${\bf u}_1=(0, 1)$, ${\bf u}_2=(1, 1)$, ${\bf u}_3=({3\over 2}, {1\over 2})$, ${\bf u}_4=({3\over 2}, 0)$, ${\bf u}_5=( 1,-{1\over 2})$, ${\bf u}_6=-{\bf u}_1$, ${\bf u}_7=-{\bf u}_2$, ${\bf u}_8=-{\bf u}_3$, ${\bf u}_9=-{\bf u}_4$ and ${\bf u}_{10}=-{\bf u}_5$ as the middle points of its edges.}

\medskip\noindent
{\bf Proof.} Let ${\bf v}_1$, ${\bf v}_2$, $\ldots $, ${\bf v}_{10}$ be the ten vertices of $P_{10}$ enumerated in the clock order, let $G_i$ denote the edge of $P_{10}$ with ends ${\bf v}_i$ and ${\bf v}_{i+1}$, where ${\bf v}_{11}={\bf v}_1$, and let ${\bf u}_i$ denote the middle point of $G_i$. For convenience, we write ${\bf v}_i=(x_i, y_i)$ and ${\bf u}_i=(x'_i, y'_i).$

It is known that $\sigma (D)+\sigma (\Lambda )$ is a $k$-fold lattice tiling of $\mathbb{E}^2$ whenever $D+\Lambda $ is such a tiling and $\sigma $ is a non-singular linear transformation from $\mathbb{E}^2$ to $\mathbb{E}^2$. Therefore, without loss of generality, we assume that $\Lambda =\mathbb{Z}^2$ and $P_{10}+\Lambda $ is a five-fold lattice tiling of $\mathbb{E}^2$.

By Lemma 1 we know that
$${\rm int}(G_i)\cap \mbox{${1\over 2}$}\Lambda \not=\emptyset $$
holds for all the ten edges $G_i$ and, if ${\bf u}_i\not \in {1\over 2}\Lambda $, then $G_i$ is a lattice vector of $\Lambda $. Now, we consider two cases.

\medskip\noindent
{\bf Case 1.} {\it $G_1$ is a lattice vector of $\Lambda $.} Without loss of generality, by a uni-modular linear transformation, we assume that ${\bf v}_2-{\bf v}_1=(k, 0)$ and $y'_1>0$, where $k$ is a positive integer. In fact, by reduction, one may assume that $G_1$ is primitive as a lattice vector and therefore $k=1$. Then, it can be deduced that
$$y_1=y'_1=y_2\in \mbox{${1\over 2}$}\mathbb{Z}$$
and all $y_i-y_{i+1}$ are integers. In particular, when $i=2$, $3,$ $4$ and $5$, they are positive integers. Thus, one can deduce that
$$y'_1\ge 2.\eqno(8)$$

\noindent
{\bf Case 1.1.} $y'_1\ge 5/2$. Let $P$ denote the parallelogram with vertices ${\bf v}_1$, ${\bf v}_2$, ${\bf v}_6$ and ${\bf v}_7$, one can deduce that
$${\rm vol}(P_{10})> {\rm vol}(P)\ge 5$$
and therefore, since $\tau^*(P_{10})$ is an integer,
$$\tau^*(P_{10})={\rm vol}(P_{10})\ge 6,\eqno(9)$$
which contradicts the assumption that $P_{10}+\Lambda $ is a five-fold tiling of $\mathbb{E}^2$.

\medskip\noindent
{\bf Case 1.2.} $y'_1=2$. Then we must have
$$y_2-y_3=y_3-y_4=y_4-y_5=y_5-y_6=1.$$
By the second term of Lemma 1, one can deduce that
$${\bf u}_i\in \mbox{${1\over 2}$} \Lambda,\quad i=2,\ 3,\ 4,\ 5. $$
Since ${\bf v}_2=(1,0)+{\bf v}_1$ and
$${\bf v}_{i+1}=2{\bf u}_i-{\bf v}_i$$
holds for all $i=2,$ $3,$ $4$ and $5$, it can be deduced that
$$-{\bf v}_1={\bf v}_6=2({\bf u}_5-{\bf u}_4+{\bf u}_3-{\bf u}_2)+(1,0)+{\bf v}_1$$
and therefore
$${\bf v}_i\in \mbox{${1\over 2}$} \Lambda, \quad i=1,\ 2,\ \ldots,\ 10. $$

\begin{figure}[!ht]
\centering
\includegraphics[scale=0.5]{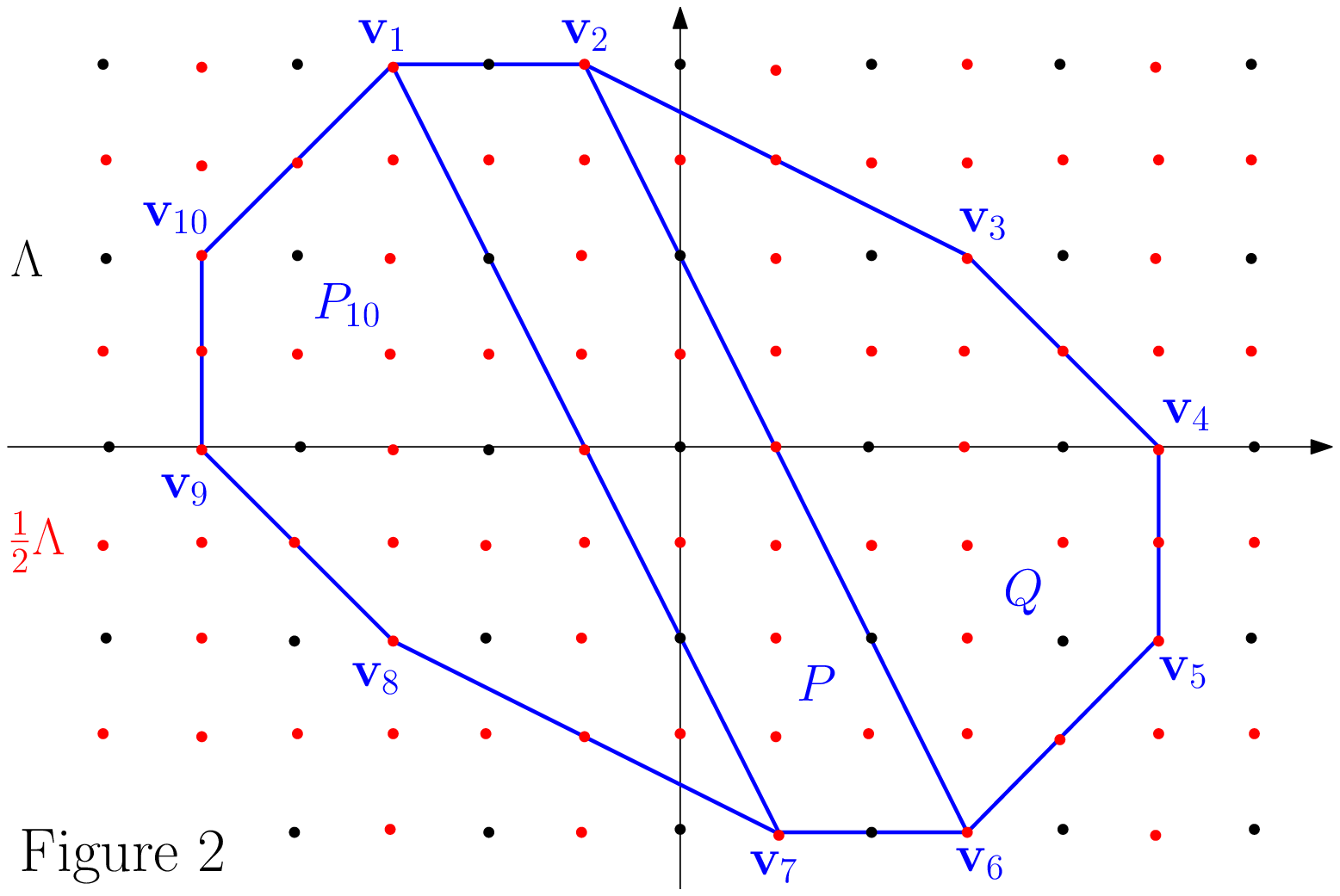}
\end{figure}

Let $P$ denote the parallelogram with vertices ${\bf v}_1$, ${\bf v}_2$, ${\bf v}_6$ and ${\bf v}_7$, and let $Q$ denote the pentagon with vertices ${\bf v}_2$, ${\bf v}_3$, ${\bf v}_4$, ${\bf v}_5$ and ${\bf v}_6$, as shown by Figure 2. Applying Pick's theorem to $Q$ and ${1\over 2}\Lambda$, we get
$${\rm vol}(Q)\ge {1\over 4}\left({9\over 2}-1\right),$$
$${\rm vol}(P_{10})={\rm vol}(P)+2\cdot {\rm vol}(Q)\ge 4+{2\over 4}\cdot \left( {9\over 2}-1\right)$$
and therfore
$$\tau^*(P_{10})={\rm vol}(P_{10})\ge 6,\eqno(10)$$
which contradicts the assumption that $P_{10}+\Lambda $ is a five-fold tiling of $\mathbb{E}^2$.

\medskip\noindent
{\bf Case 2.} {\it All the middle points ${\bf u}_i$ belong to ${1\over 2}\Lambda $.} Since $P_{10}+\Lambda $ is a five-fold lattice tiling of $\mathbb{E}^2$, one can deduce that
$${\rm vol}(2P_{10})=20$$
and all ${\bf u}'_i=2{\bf u}_i$ belong to $\Lambda $. For convenience, we define $Q_{10}$ to be the centrally symmetric lattice decagon with vertices ${\bf u}'_1$, ${\bf u}'_2$, $\ldots $, ${\bf u}'_{10}$ and write ${\bf u}'_i=(x'_i, y'_i)$. Since $Q_{10}$ is a centrally symmetric lattice polygon, its area must be a positive integer. Thus, we have
$${\rm vol}(Q_{10})\le 19.\eqno(11)$$

\medskip
Now, we explore $Q_{10}$ in detail by considering the following subcases.

\medskip\noindent
{\bf Case 2.1.} {\it ${\bf u}'_1$ is primitive in $\Lambda $}. Without loss of generality, guaranteed by uni-modular linear transformations, we take ${\bf u}'_1=(0, 1)$.  Then, Lemma 3 implies
$$\left\{
\begin{array}{ll}
x'_4-x'_5&\hspace{-0.3cm}=x'_3-x'_2, \\
y'_4-y'_5&\hspace{-0.3cm}=y'_3-y'_2+1.
\end{array}\right. \eqno(12) $$

If $x'_2\ge x'_3$ or $x'_3=x'_4$, one can easily deduce contradiction with convexity from (12). For example, if $x'_3=x'_4>x'_2$, then it can be deduced by (12) that
$${\bf u}'_2-{\bf u}'_5={\bf u}'_{10}-{\bf u}'_7=k{\bf u}'_1$$
with $k\ge 2$, which contradicts the assumption that $Q_{10}$ is a centrally symmetric convex decagon. Therefore, we may assume that
$$x'_3> x'_i\eqno(13)$$
for all $i\not= 3.$

Let $T'$ denote the lattice triangle with vertices ${\bf u}'_1$, ${\bf u}'_2$ and ${\bf u}'_3$, let $Q$ denote the lattice quadrilateral with vertices ${\bf u}'_3$, ${\bf u}'_4$, ${\bf u}'_5$ and ${\bf u}'_6$, and let $T$ denote the lattice triangle with vertices ${\bf u}'_1$, ${\bf u}'_3$ and ${\bf u}'_6$ (as shown by Figure 3). It follows from (11) and Pick's theorem that
$${\rm vol}(T)\le {1\over 2}\Bigl( 19 -2 \bigl({\rm vol}(T')+{\rm vol}(Q)\bigr)\Bigr)\le 8$$
and therefore
$$x'_3\le 8.\eqno(14)$$

\begin{figure}[!ht]
\centering
\includegraphics[scale=0.6]{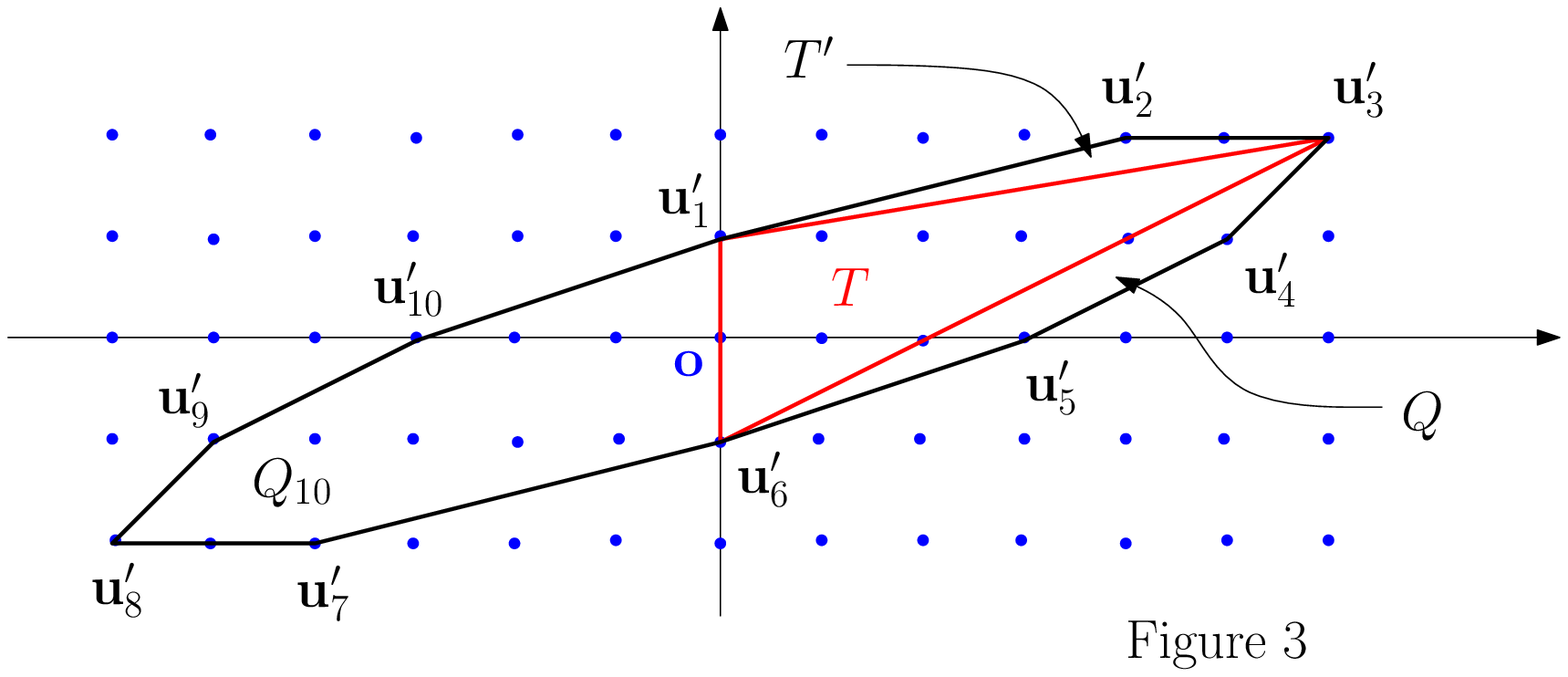}
\end{figure}

Let $\alpha $ denote the slope of $G_1$, that is
$$\alpha ={{y_2-y_1}\over {x_2-x_1}}.$$
By a uni-modular linear transformation such as
$$\left\{
\begin{array}{ll}
x'=x,&\\
y'=y+kx,&
\end{array}
\right.$$
where $k$ is a suitable integer, we may assume that
$$0\le \alpha <1.\eqno (15)$$
Let $L_i$ denote the straight line containing $G_i$, it is obvious that $P_{10}$ is in the strip bounded by $L_1$ and $L_6$. Furthermore, we define
five slopes
$$\beta_i={{y'_{i+1}-y'_i}\over {x'_{i+1}-x'_i}},\quad i=1,\ 2,\ \ldots,\ 5.$$

By convexity it can be shown that there is no five-fold lattice decagon tile with $\alpha =0$ in our setting. When $\alpha >0$, by (12) and convexity it follows that $y'_4-y'_5\ge 1$ and therefore
$$y'_3-y'_2\ge 0.\eqno (16)$$

\begin{figure}[!ht]
\centering
\includegraphics[scale=0.6]{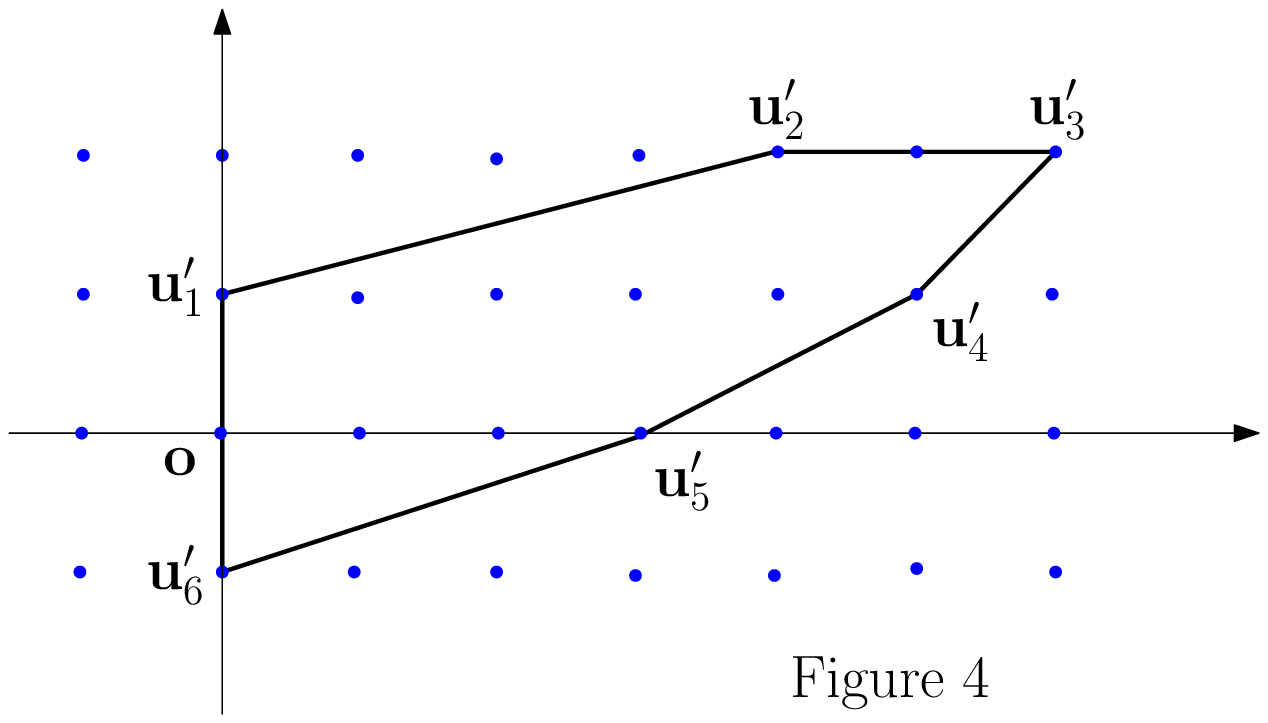}
\end{figure}

As shown by Figure 4, we assume that $${\bf u}'_3-{\bf u}'_4=(p_1, q_1)$$ and $${\bf u}'_5-{\bf u}'_6=(p_2, q_2),$$
where all $p_i$ and $q_i$ are positive integers. Then, by (14) we have
$$x'_3-x'_2=x'_3-(x'_2-x'_1)=x'_3-(p_1+p_2)\le 6.$$
Now, we consider in subcases with respect to the different orientations of ${\bf u}'_3-{\bf u}'_2.$

\medskip\noindent
{\bf Case 2.1.1.} {\it $y'_3-y'_2=0$ and $x'_3-x'_2=1$}. By (12) and convexity we have $x'_4-x'_5=1$, $y'_4-y'_5=1$, $\beta_4=1$,
$$\beta_3={{q_1}\over {p_1}}>1$$
and
$$\beta_5={{q_2}\over {p_2}}<1.$$
Then, one can deduce that
$$\beta_1={{q_1+q_2-1}\over {p_1+p_2}}>{{q_2}\over {p_2}}=\beta_5,$$
which contradicts the convexity of $Q_{10}$.

\medskip\noindent
{\bf Case 2.1.2.} {\it $y'_3-y'_2=0$ and $x'_3-x'_2=2$}. By (12) and convexity we have $x'_4-x'_5=2$, $y'_4-y'_5=1$, $\beta_4={1\over 2}$,
$$\beta_3={{q_1}\over {p_1}}>{1\over 2},\eqno(17)$$
$$\beta_5={{q_2}\over {p_2}}<{1\over 2}\eqno(18)$$
and
$$\beta_1={{q_1+q_2-1}\over {p_1+p_2}}<{{q_2}\over {p_2}}.\eqno(19)$$

By (14) and (18) one can deduce that
$$3\le p_2\le 5,\eqno(20)$$
$$1\le p_1\le 3\eqno(21)$$
and
$$1\le q_2\le 2.\eqno(22)$$
On the other hand, by (19), (21) and (22) we get
$$p_2(q_1-1)<p_1q_2\le 6$$
and therefore
$$1\le q_1\le 2.\eqno(23)$$

Then, it can be verified that the only integer groups $(p_1, q_1, p_2,q_2)$ satisfying (14), (17), (18) and (19) are $(1,1,3,1),$ $(1,1,4,1)$, $(1,1,5,1)$ and $(1,1,5,2)$. By checking the areas of their corresponding decagons, keeping the subcase conditions in mind, the only $Q_{10}$ satisfying (11) is the one with vertices ${\bf u}'_1=(0,1)$, ${\bf u}'_2=(4,2)$, ${\bf u}'_3=(6,2)$, ${\bf u}'_4=(5,1)$, ${\bf u}'_5=(3,0)$, ${\bf u}'_6=-{\bf u}'_1$, ${\bf u}'_7=-{\bf u}'_2$, ${\bf u}'_8=-{\bf u}'_3$, ${\bf u}'_9=-{\bf u}'_4$ and ${\bf u}'_{10}=-{\bf u}'_5$. Clearly, this lattice polygon is equivalent to the one stated in the theorems under the linear transformation
$$\left\{
\begin{array}{ll}
x'={1\over 2}(x-2y),&\\
\vspace{-0.3cm}
&\\
y'={1\over 2}y.&
\end{array}
\right.$$

\medskip\noindent
{\bf Case 2.1.3.} {\it $y'_3-y'_2=0$ and $x'_3-x'_2=3$}. By (12) and convexity we have $x'_4-x'_5=3$, $y'_4-y'_5=1$, $\beta_4={1\over 3}$,
$$\beta_3={{q_1}\over {p_1}}>{1\over 3},$$
$$\beta_5={{q_2}\over {p_2}}<{1\over 3}\eqno(24)$$
and
$$\beta_1={{q_1+q_2-1}\over {p_1+p_2}}<{{q_2}\over {p_2}}.\eqno(25)$$

By (14), (24) and (25), it can be deduced that $p_2=4$, $q_2=1$, $p_1=1$ and $q_1=1$. Then we have
$${\rm vol}(Q_{10})=25,\eqno(26)$$
which contradicts (11).

\medskip\noindent
{\bf Case 2.1.4.} {\it $y'_3-y'_2=0$ and $x'_3-x'_2\ge 4$}. Then, one can easily deduce that $\beta_5<{1\over 4}$, $p_2>4$ and
$$p_1+p_2+4>8,$$
which contradicts the restriction of (14).

\medskip\noindent
{\bf Case 2.1.5.} {\it $y'_3-y'_2=1$ and $x'_3-x'_2=1$}. Then, by convexity we get
$$\alpha >\beta_1>\beta_2=1,$$
which contradicts the assumption of (15).

\medskip\noindent
{\bf Case 2.1.6.} {\it $y'_3-y'_2=1$ and $x'_3-x'_2=2$}. By (12) and convexity we get $x'_4-x'_5=2$, $y'_4-y'_5=2,$ $\beta_4=1$,
$$\beta_3={{q_1}\over {p_1}}>1$$
and
$$\beta_5={{q_2}\over {p_2}}<1.$$
Then, it can be deduced that
$$\beta_1={{q_1+q_2-1}\over {p_1+p_2}}> {{q_2}\over {p_2}}=\beta_5,$$
which contradicts the convexity assumption of $Q_{10}$.

\medskip\noindent
{\bf Case 2.1.7.} {\it $y'_3-y'_2=1$ and $x'_3-x'_2=3$}. Then we have $x'_4-x'_5=3$, $y'_4-y'_5=2,$ $\beta_2={1\over 3}$ and $\beta_4={2\over 3}$.

On one hand, by (14) it follows that $p_2\le 4$. On the other hand, by $\beta_2<\beta_1<\beta_5<\beta_4$ it follows that
$${1\over 3}<{{q_2}\over {p_2}}<{2\over 3}.$$
Thus, the integer pair $(p_2, q_2)$ has only two choices $(2,1)$ and $(4,2)$.

Then, by checking
$${{q_1}\over {p_1}}>{2\over 3},$$
$${1\over 3}<{{q_1+q_2-1}\over {p_1+p_2}}<{1\over 2}$$
and
$$p_1+p_2\le 5,$$
it can be deduced that the only candidate for $(p_1,q_1, p_2,q_2)$ is $(1,1,4,2)$. Unfortunately, in this case we have
$${\rm vol}(Q_{10})=22,\eqno(27)$$
which contradicts the restriction of (11).

\medskip\noindent
{\bf Case 2.1.8.} {\it $y'_3-y'_2=1$ and $x'_3-x'_2=4$}. By (12), (14) and convexity it can be deduced that $p_2\le 3$, $\beta_4={1\over 2}$ and $\beta_5<\beta_4.$ Consequently, we have $p_2=3$, $q_2=1$ and $\beta_5={1\over 3}.$ Thus, by $\beta_2={1\over 4}$ and $\beta_2<\beta_1<\beta_5$ we get
$${1\over 4}<{{q_1+q_2-1}\over {p_1+p_2}}<{1\over 3}.\eqno (28)$$
However, by (14) we have $p_1+p_2\le 4$ and therefore (28) has no integer solution.

\medskip\noindent
{\bf Case 2.1.9.} {\it $y'_3-y'_2=1$ and $x'_3-x'_2\ge 5$}. It follows by (14) that $p_2\le 2$. Then we get both $\beta_4\le {2\over 5}$ and $\beta_5\ge {1\over 2}$, which contradicts the convexity of $Q_{10}$.

\medskip\noindent
{\bf Case 2.1.10.} {\it $y'_3-y'_2=2$ and $x'_3-x'_2=3$}. Then by (12) and convexity we get $\beta_2={2\over 3}$, $\beta_4=1$, $\beta_2<\beta_5<\beta_4$ and therefore
$${2\over 3}<{{q_2}\over {p_2}}<1.\eqno(29)$$
Clearly, by (14) we have $p_2\le 4$ and therefore (29) has only one group of integer solutions $p_2=4$ and $q_2=3$. Then, $\beta_2<\beta_1<\beta_5$ can be reformulated as
$${2\over 3}<{{q_1+2}\over 5}<{3\over 4},$$
which has no integer solution.

\medskip\noindent
{\bf Case 2.1.11.} {\it $y'_3-y'_2=2$ and $x'_3-x'_2=4$}. Then by (12) and convexity we get $\beta_2={1\over 2}$, $\beta_4={3\over 4}$, $\beta_2<\beta_5<\beta_4$ and therefore
$${1\over 2}<{{q_2}\over {p_2}}<{3\over 4}.\eqno(30)$$
Clearly, by (14) we have $p_2\le 3$ and therefore (30) has only one group of integer solutions $p_2=3$ and $q_2=2$. Then, $\beta_2<\beta_1<\beta_5$ can be reformulated as
$${1\over 2}<{{q_1+1}\over 4}<{2\over 3},$$
which has no integer solution.

\medskip\noindent
{\bf Case 2.1.12.} {\it $y'_3-y'_2=2$ and $x'_3-x'_2=5$}. Then by (12) and convexity we get $\beta_2={2\over 5}$, $\beta_4={3\over 5}$, $\beta_2<\beta_5<\beta_4$ and therefore
$${2\over 5}<{{q_2}\over {p_2}}<{3\over 5}.\eqno(31)$$
Clearly, by (14) we have $p_2\le 2$ and therefore (31) has only one group of integer solutions $p_2=2$ and $q_2=1$. Then, $\beta_2<\beta_1<\beta_5$ can be reformulated as
$${2\over 5}<{{q_1}\over 3}<{1\over 2},$$
which has no integer solution.

\medskip\noindent
{\bf Case 2.1.13.} {\it $y'_3-y'_2=2$ and $x'_3-x'_2=6$}. Then by (12) and convexity we get $\beta_4={1\over 2}$ and $\beta_5\ge 1$, which contradicts the convexity of $Q_{10}$.

\medskip\noindent
{\bf Case 2.1.14.} {\it $y'_3-y'_2=3$ and $x'_3-x'_2=4$}. Then by (12) and convexity we have $p_2\le 3$, $\beta_2={3\over 4},$ $\beta_4=1$ and $\beta_2<\beta_5<\beta_4$. Unfortunately, the inequalities
$${3\over 4}<{{q_2}\over {p_2}}<1$$
and
$$p_2\le 3$$
have no common integer solution.

\medskip\noindent
{\bf Case 2.1.15.} {\it $y'_3-y'_2=3$ and $x'_3-x'_2=5$}. Then by (12) and convexity we have $p_2\le 2$, $\beta_2={3\over 5},$ $\beta_4={4\over 5}$ and $\beta_2<\beta_5<\beta_4$. However, the inequalities
$${3\over 5}<{{q_2}\over {p_2}}<{4\over 5}$$
and
$$p_2\le 2$$
have no common integer solution.

\medskip\noindent
{\bf Case 2.1.16.} {\it $y'_3-y'_2=3$ and $x'_3-x'_2=6$}. Then by (12) and convexity we get $\beta_4={2\over 3}$ and $\beta_5\ge 1$, which contradicts the convexity of $Q_{10}$.

\medskip\noindent
{\bf Case 2.1.17.} {\it $y'_3-y'_2=4$ and $x'_3-x'_2=5$}. Then by (12) and convexity we have $p_2\le 2$, $\beta_2={4\over 5},$ $\beta_4=1$ and $\beta_2<\beta_5<\beta_4$. Unfortunately, the inequalities
$${4\over 5}<{{q_2}\over {p_2}}<1$$
and
$$p_2\le 2$$
have no common integer solution.

\medskip\noindent
{\bf Case 2.1.18.} {\it $y'_3-y'_2=4$ and $x'_3-x'_2=6$}. Then by (12) and convexity we get $\beta_4={5\over 6}$ and $\beta_5\ge 1$, which contradicts the convexity of $Q_{10}$.

\medskip\noindent
{\bf Case 2.1.19.} {\it $y'_3-y'_2=5$ and $x'_3-x'_2=6$}. Then by (12) and convexity we get $\beta_4=1$ and $\beta_5\ge 1$, which contradicts the convexity of $Q_{10}$.

\medskip\noindent
{\bf Case 2.2.} {\it All ${\bf u}'_i$ are even multiplicative}. Then all ${\bf u}_i$ belong to $\Lambda $. It follows by Lemma 1 that ${1\over 2}P_{10}+\Lambda $ is a $k$-fold lattice tiling with
$$k={\rm vol}\left(\mbox{${1\over 2}$}P_{10}\right)={5\over 4},$$
which contradicts the assumption that $k$ is a positive integer.

\medskip\noindent
{\bf Case 2.3.} {\it All ${\bf u}'_i$ are multiplicative, ${\bf u}'_1$ is odd multiplicative}. Without loss of generality, guaranteed by uni-modular linear transformations, we take ${\bf u}'_1=(0, 2q+1)$, where $q$ is a positive integer.

By Lemma 3 it follows that
$$x'_4-x'_5=x'_3-x'_2.$$
Therefore, by convexity and reflection we may assume that
$$x'_3\ge x'_i, \qquad i=1, 2, \ldots , 10.$$

Let $T'$ denote the lattice triangle with vertices ${\bf u}'_1$, ${\bf u}'_2$ and ${\bf u}'_3$, let $Q$ denote the lattice quadrilateral with vertices ${\bf u}'_3$, ${\bf u}'_4$, ${\bf u}'_5$ and ${\bf u}'_6$, and let $T$ denote the lattice triangle with vertices ${\bf u}'_1$, ${\bf u}'_3$ and ${\bf u}'_6$, as shown in Figure 3. It follows from (11) and Pick's theorem that
$${\rm vol}(T)\le {1\over 2}\Bigl( 19 -2 \bigl({\rm vol}(T')+{\rm vol}(Q)\bigr)\Bigr)\le 8$$
and therefore
$$x'_3 ={{2\cdot {\rm vol}(T)}\over {2(2q+1)}}\le \left\lfloor {8\over 3}\right\rfloor = 2.$$

It is assumed that all ${\bf u}'_i$ are multiplicative. Therefore we have
$$x'_2=x'_3=x'_4=x'_5=2,$$
which contradicts the convexity of $Q_{10}$.

\medskip
As a conclusion of all these cases, Lemma 8 is proved.\hfill{$\Box$}

\medskip\noindent
{\bf Lemma 9.} {\it For every centrally symmetric convex octagon $P_8$ we have
$$\tau^*(P_8)\ge 5,$$
where the equality holds if and only if, under a suitable affine linear transformation, it with vertices ${\bf v}_1=(-\alpha , -{3\over 2})$, ${\bf v}_2=(1-\alpha , -{3\over 2})$, ${\bf v}_3=(1+\alpha , -{1\over 2})$, ${\bf v}_4=(1-\alpha , {1\over 2})$, ${\bf v}_5=-{\bf v}_1$, ${\bf v}_6=-{\bf v}_2$, ${\bf v}_7=-{\bf v}_3$ and ${\bf v}_8=-{\bf v}_4$, where $0<\alpha <{1\over 4}$, or with vertices ${\bf v}_1=(\beta , -2)$, ${\bf v}_2=(1+\beta , -2)$, ${\bf v}_3=(1-\beta , 0)$, ${\bf v}_4=(\beta , 1)$, ${\bf v}_5=-{\bf v}_1$, ${\bf v}_6=-{\bf v}_2$, ${\bf v}_7=-{\bf v}_3$, ${\bf v}_8=-{\bf v}_4$, where ${1\over 4}<\beta <{1\over 3}$.}

\medskip\noindent
{\bf Proof.} Let $P_8$ be a centrally symmetric convex octagon centered at the origin, let ${\bf v}_1$, ${\bf v}_2$, $\ldots$, ${\bf v}_8$ be the eight vertices of $P_8$ enumerated in an anti-clock order, let $G_i$ denote the edge with ends ${\bf v}_i$ and ${\bf v}_{i+1}$, where ${\bf v}_9={\bf v}_1$, and let ${\bf u}_i$ denote the midpoint of $G_i$. For convenience, we write ${\bf v}_i=(x_i,y_i)$ and ${\bf u}_i=(x'_i,y'_i)$. Assume that $\Lambda = \mathbb{Z}^2$ and $P_8+\Lambda $ is a five-fold lattice tiling. Then, we have
$$\tau^*(P_8)={\rm vol}(P_8)=5.\eqno(32)$$

Based on Lemma 2, by a uni-modular transformation, we may assume that $G_1\cap {1\over 2}\Lambda\not=\emptyset$ and ${\bf v}_2-{\bf v}_1=(k, 0),$
where $k$ is a positive integer.  If $k>1$, we define $P_8'$ to be the octagon with vertices ${\bf v}'_1={\bf v}_1+({{k-1}\over 2},0),$ ${\bf v}'_2={\bf v}_2+({{1-k}\over 2},0),$ ${\bf v}'_3={\bf v}_3+({{1-k}\over 2},0),$ ${\bf v}'_4={\bf v}_4+({{1-k}\over 2},0),$ ${\bf v}'_5={\bf v}_5+({{1-k}\over 2},0),$ ${\bf v}'_6={\bf v}_6+({{k-1}\over 2},0),$ ${\bf v}'_7={\bf v}_7+({{k-1}\over 2},0)$ and ${\bf v}'_8={\bf v}_8+({{k-1}\over 2},0)$, as shown by Figure 5. By Lemma 1 it can be shown that
$P'_8+\Lambda $ is a multiple lattice tiling of $\mathbb{E}^2$ and therefore
$$\tau^*(P'_8)\le {\rm vol}(P'_8)<{\rm vol}(P_8)=5,$$
which contradicts the known fact that $\tau^* (P'_8)\ge 5$. Thus, we have ${\bf v}_2-{\bf v}_1=(1,0).$

\begin{figure}[!ht]
\centering
\includegraphics[scale=0.52]{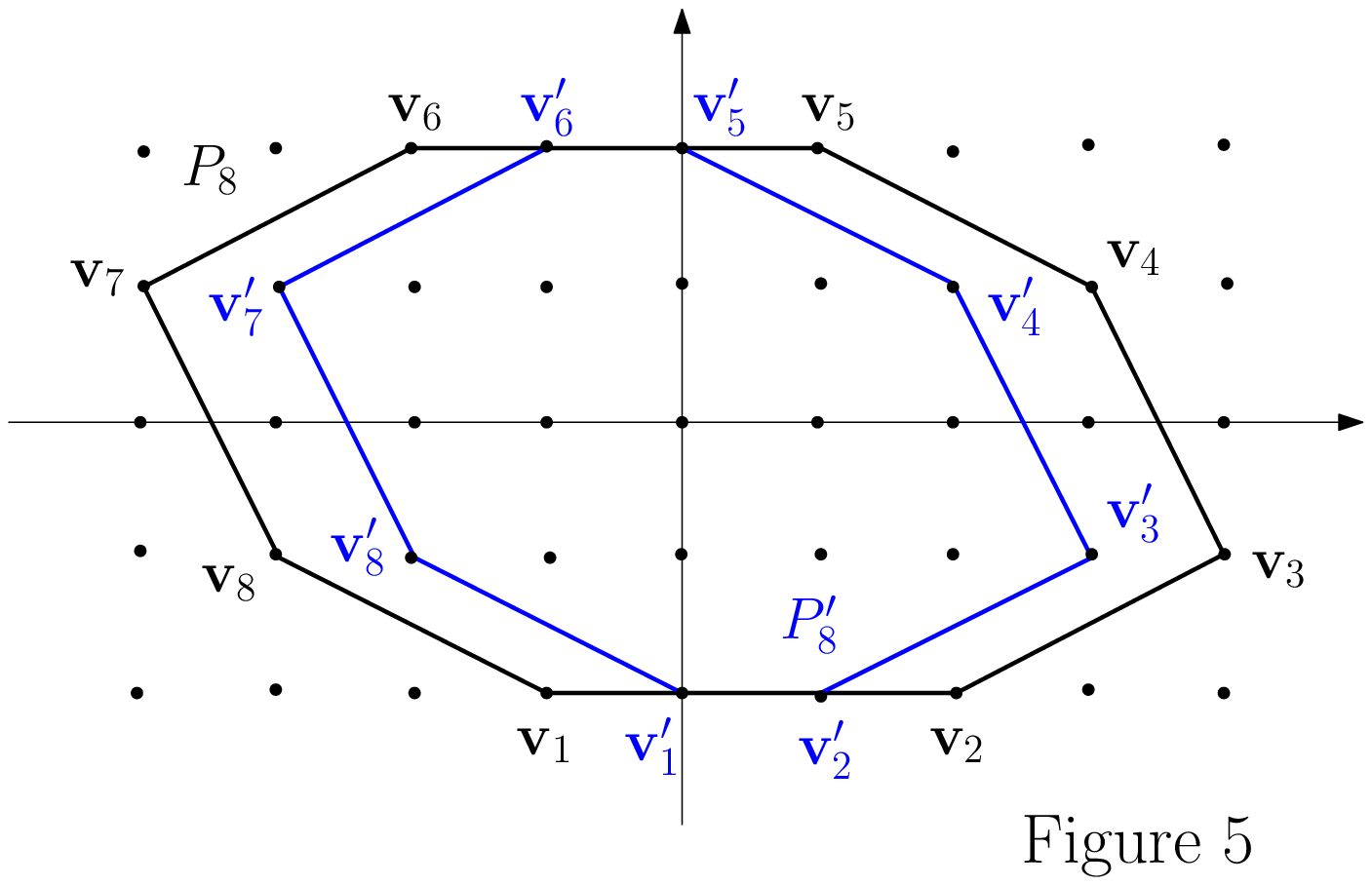}
\end{figure}

Apply Lemma 1 successively to $G_1$, $G_2$, $G_3$ and $G_4$, one can deduce that all $2y_2$, $y_3-y_2$, $y_4-y_3$ and $y_5-y_4$ are positive integers. Therefore, we have
$$y_2=y_1\le -{3\over 2}.$$
On the other hand, if $y_2=y_1\le -{5\over 2}$ and let $P$ denote the parallelogram with vertices ${\bf v}_1$, ${\bf v}_2$, ${\bf v}_5$ and ${\bf v}_6$, it can be deduced by convexity that
$${\rm vol}(P_8)>{\rm vol}(P)\ge 5,$$
which contradicts the assumption of (32). Thus, to prove the theorem it is sufficient to deal with the two cases
$$y_2=y_1=-{3\over 2}, \ -2.$$

\noindent
{\bf Case 1.} $y_2=y_1=-{3\over 2}.$  In this case,
$$y_{i+1}-y_i=1$$
must hold for all $i=2, 3$ and $4$. Then, it follows by Lemma 1 that all the midpoints of $G_2$, $G_3$ and $G_4$ belong to ${1\over 2}\Lambda $. Furthermore, by a uni-modular transformation
$$\left\{
\begin{array}{ll}
x'\hspace{-0.2cm}&=x-ky,\\
y'\hspace{-0.2cm}&=y,
\end{array}\right.$$
with a suitable integer $k$, we may assume that $-{5\over 4}\le x_1<{1\over 4}$.

If $G_2$ is vertical, then $x_2$ is an integer or an half integer. Consequently, we have $x_1\in {1\over 2}\mathbb{Z}$. Therefore $x_1$ only can be $-1$, $-{1\over 2}$ or $0$. By considering three subcases with respect to $x_1=-1$, $-{1\over 2}$ or $0$, it can be deduced that there is no octagon of this type satisfying Lemma 1. For example, when $x_1=-{1\over 2}$, by Lemma 1 and convexity we have ${\bf v}_1=\left(-{1\over 2}, -{3\over 2}\right)$, ${\bf v}_2=\left({1\over 2}, -{3\over 2}\right)$, ${\bf v}_3=\left({1\over 2}, -{1\over 2}\right)$, ${\bf v}_4=\left({1\over 2}, {1\over 2}\right)$, ${\bf v}_5=-{\bf v}_1$, ${\bf v}_6=-{\bf v}_2$, ${\bf v}_7=-{\bf v}_3$ and ${\bf v}_8=-{\bf v}_4$. Then, $P_8$ is no longer an octagon but a parallelogram.

\begin{figure}[!ht]
\centering
\includegraphics[scale=0.5]{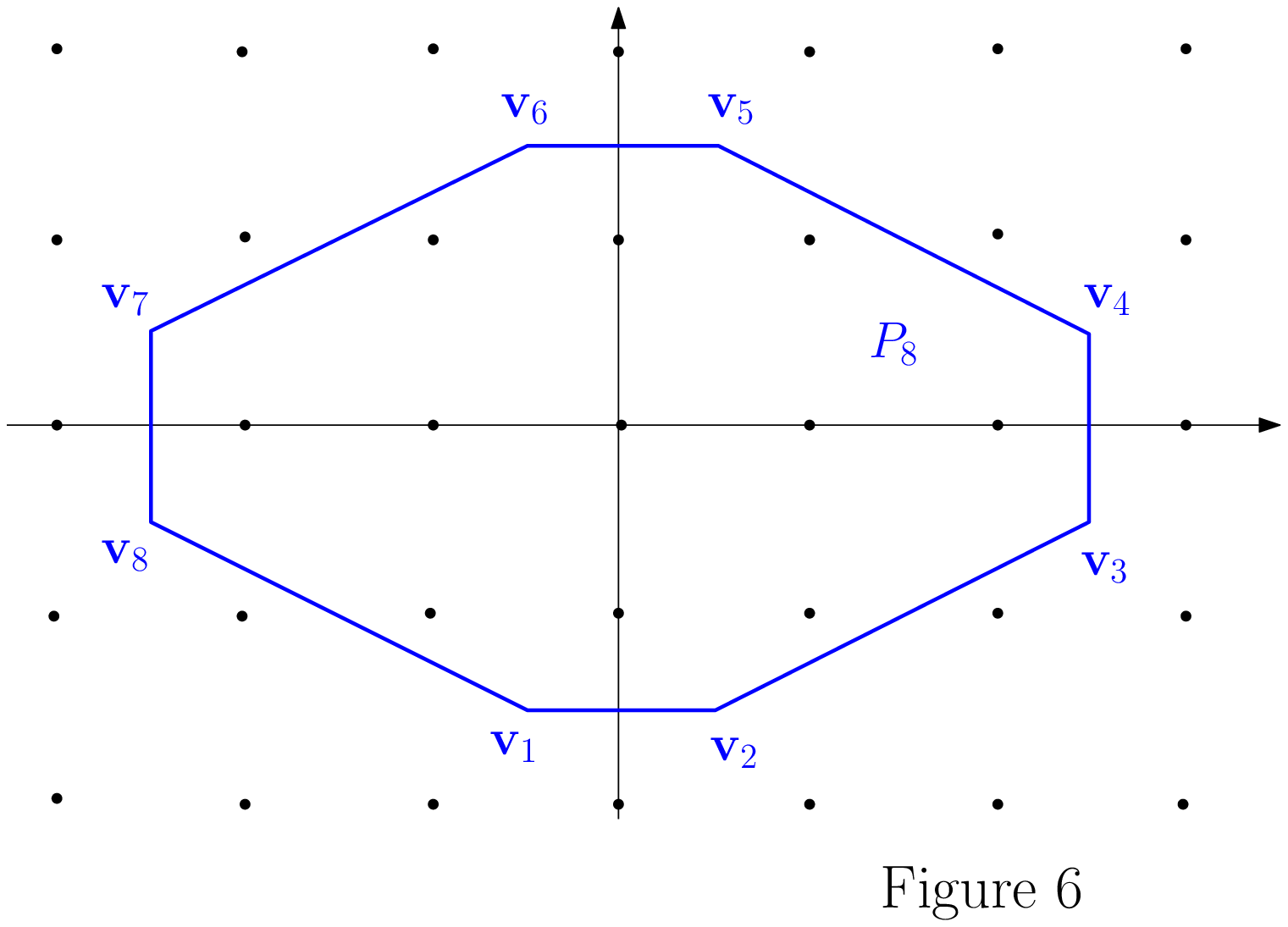}
\end{figure}

If $G_3$ is vertical, then $x_1$ must be an integer or an half integer as well. Therefore, it only can be $-1$, $-{1\over 2}$ or $0$. By considering three subcases with respect to $x_1=-1$, $-{1\over 2}$ or $0$, it can be deduced that
$${\rm vol}(P_8)\ge 7,\eqno(33)$$
which contradicts the assumption of (32). For example, when $x_1=-{1\over 2}$, by Lemma 1 and convexity we have ${\bf v}_1=\left(-{1\over 2}, -{3\over 2}\right)$, ${\bf v}_2=\left({1\over 2}, -{3\over 2}\right)$, ${\bf v}_3=\left({1\over 2}+k, -{1\over 2}\right)$, ${\bf v}_4=\left({1\over 2}+k, {1\over 2}\right)$, ${\bf v}_5=-{\bf v}_1$, ${\bf v}_6=-{\bf v}_2$, ${\bf v}_7=-{\bf v}_3$ and ${\bf v}_8=-{\bf v}_4$, where $k$ is a positive integer. Then, as shown by Figure 6, it can be deduced that
$${\rm vol}(P_8)=3+4 k\ge 7.$$

If none of the three edges $G_2$, $G_3$ and $G_4$ is vertical, by convexity it is sufficient to deal with the following three subcases.

\medskip
\noindent
{\bf Subcase 1.1.}  $x'_3>\max\{ x'_2, x'_4\}$. Then we replace the eight vertices ${\bf v}_3$, ${\bf v}_4$, ${\bf v}_5$, ${\bf v}_6$, ${\bf v}_7$, ${\bf v}_8$, ${\bf v}_1$ and ${\bf v}_2$ by ${\bf v}'_3=(x'_3, -{1\over 2})$, ${\bf v}'_4=(x'_3, {1\over 2})$,
${\bf v}'_5=(2x'_4-x'_3, {3\over 2})$, ${\bf v}'_6=(2x'_4-x'_3-1, {3\over 2})$, ${\bf v}'_7=-{\bf v}'_3$, ${\bf v}'_8=-{\bf v}'_4$,  ${\bf v}'_1=-{\bf v}'_5$ and ${\bf v}'_2=-{\bf v}'_6$, respectively (as shown by Figure 7). In practice, one first makes $G_3$ vertical and then changes the other vertices successively. Clearly, this process does not change the area of the polygon. Then one can deduce that $x'_3\ge {3\over 2}$ and therefore
$${\rm vol}(P_8)=3\cdot 2x'_3-(2x'_3-1)=4x'_3+1\ge 7,\eqno (34)$$
which contradicts the assumption of (32).

\begin{figure}[!ht]
\centering
\includegraphics[scale=0.55]{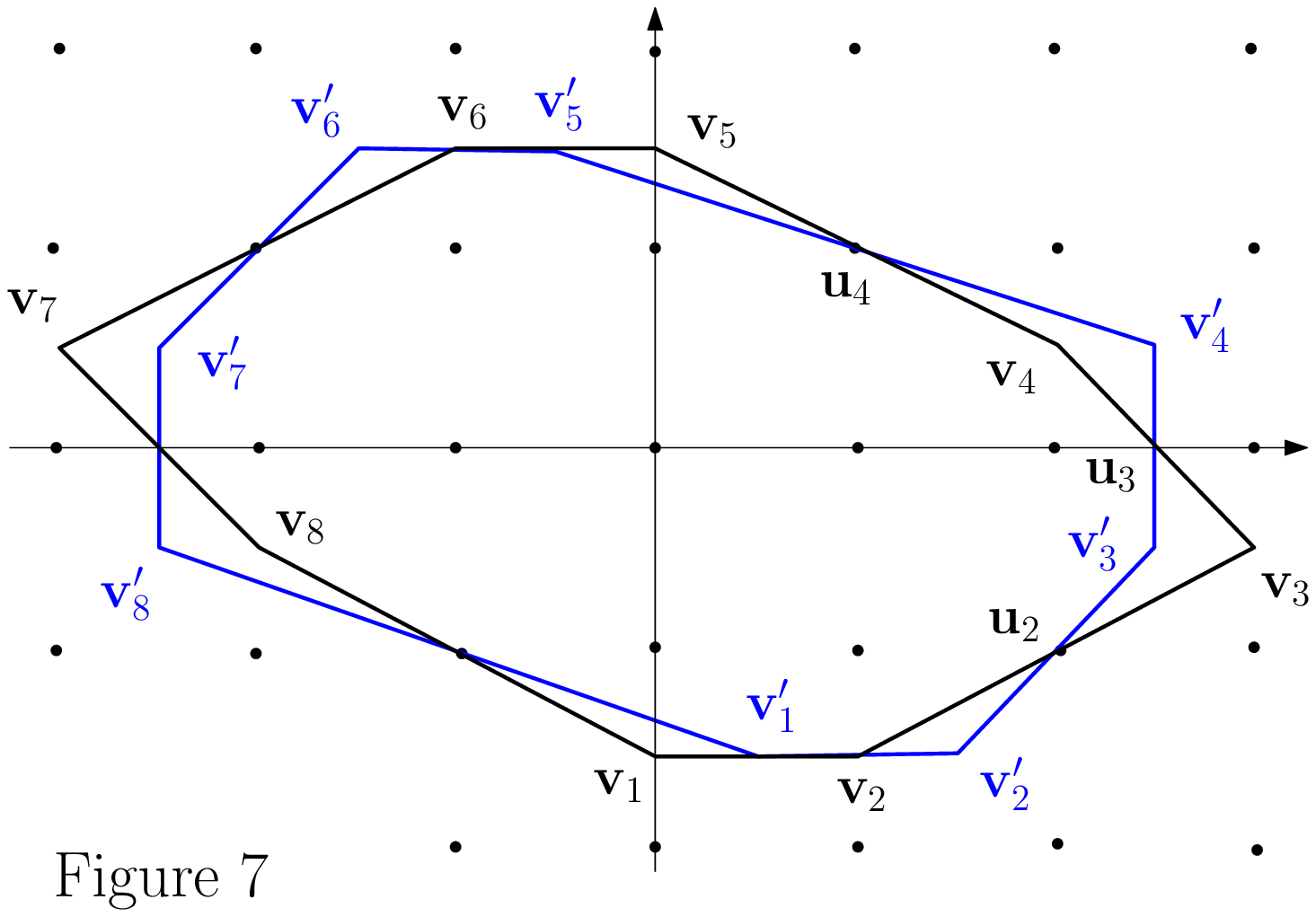}
\end{figure}

\medskip
\noindent
{\bf Subcase 1.2.} $x'_2>\max \{ x'_3, x'_4\}$. If $x_3>x_2$, one can repeat the above process. At the end we get $x'_2\ge {3\over 2}$ and
$${\rm vol}(P_8)> 3\cdot 2x'_2-2(2x'_2-1)=2x'_2+2\ge 5,\eqno(35)$$
which contradicts the assumption of (32). If $x_2>x_3$, since $-{5\over 4}\le x_1<{1\over 4}$, ${\bf u}_2$ only can be $(1,-1)$, $({1\over 2}, -1)$, $(0,-1)$ or $(-{1\over 2},-1)$. Then it can be easily checked that there is no convex octagon of this type satisfying Lemma 1.

\medskip
\noindent
{\bf Subcase 1.3.} $x'_2=x'_3>x'_4$. Then, we replace the eight vertices ${\bf v}_2$, ${\bf v}_3$, ${\bf v}_4$, ${\bf v}_5$, ${\bf v}_6$, ${\bf v}_7$, ${\bf v}_8$ and  ${\bf v}_1$ by ${\bf v}'_2=(x'_2, -{3\over 2})$, ${\bf v}'_3=(x'_2, -{1\over 2})$, ${\bf v}'_4=(x'_2, {1\over 2})$,
${\bf v}'_5=2{\bf u}_4-{\bf v}'_4$, ${\bf v}'_6=-{\bf v}'_2$, ${\bf v}'_7=-{\bf v}'_3$, ${\bf v}'_8=-{\bf v}'_4$ and  ${\bf v}'_1=-{\bf v}'_5$, respectively (as shown by Figure 8). In practice, one first makes $G_2$ and $G_3$ vertical and then changes the other vertices successively, keeping the rules of Lemma 1. Clearly, this process does not change the area of the polygon, $x'_2\ge 1$ and therefore
$${\rm vol}(P_8)=3\cdot 2x'_2-(2x'_2-1)=4x'_2+1\ge 5,\eqno (36)$$
where the equality holds if and only if $P_8$ with vertices ${\bf v}_1=(-\alpha , -{3\over 2})$, ${\bf v}_2=(1-\alpha , -{3\over 2})$,
${\bf v}_3=(1+\alpha , -{1\over 2})$, ${\bf v}_4=(1-\alpha , {1\over 2})$, ${\bf v}_5=-{\bf v}_1$, ${\bf v}_6=-{\bf v}_2$, ${\bf v}_7=-{\bf v}_3$
and ${\bf v}_8=-{\bf v}_4$, where $0<\alpha <{1\over 4}.$ They are the octagons of first type listed in Lemma 9.

\begin{figure}[!ht]
\centering
\includegraphics[scale=0.55]{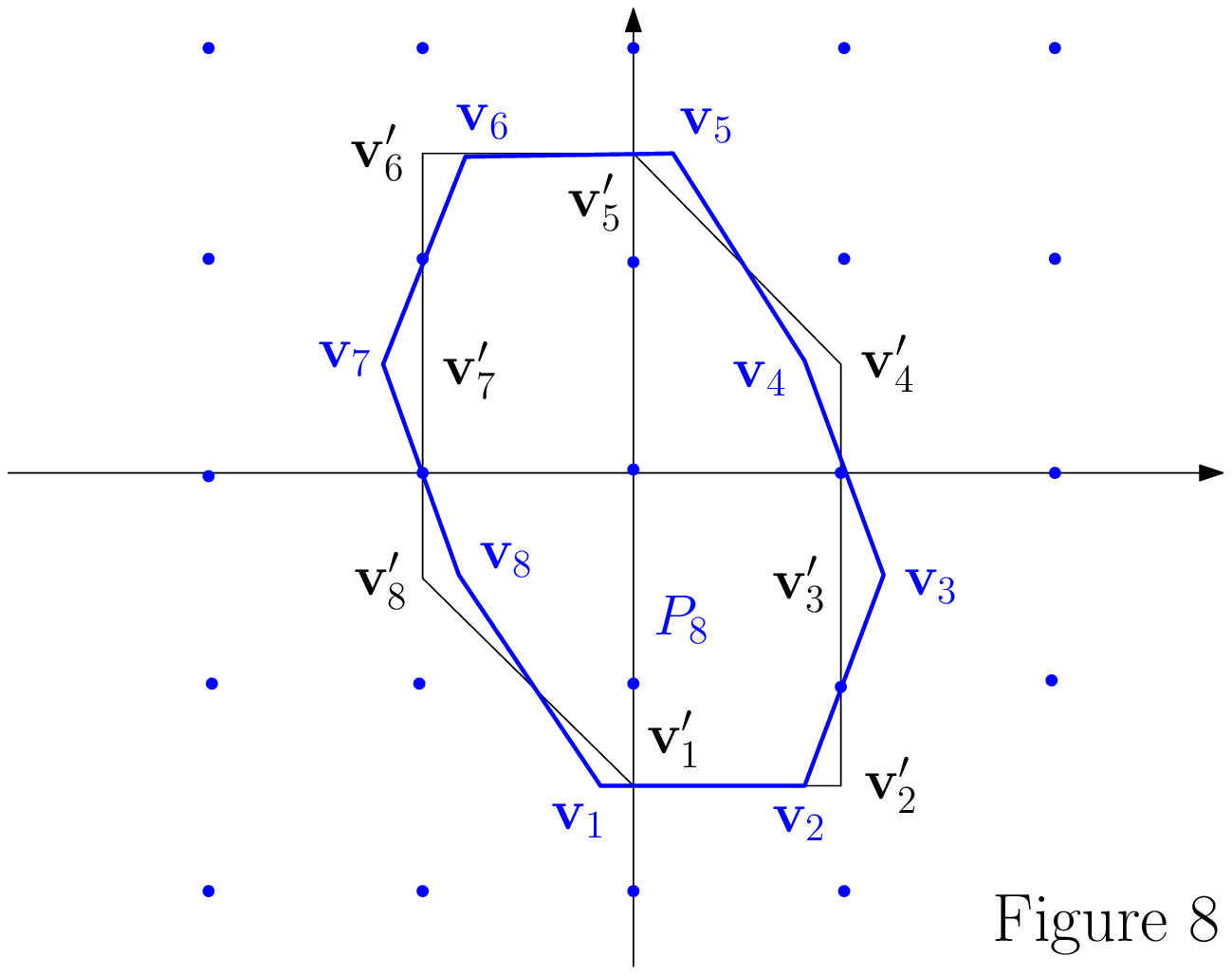}
\end{figure}

\medskip\noindent
{\bf Case 2.} $y_2=y_1=-2.$  Then, it can be deduced that one of $y_3-y_2$, $y_4-y_3$ and $y_5-y_4$ is two and the others are ones, and the midpoint ${\bf u}_i$ must belong to ${1\over 2}\Lambda $ whenever $y_{i+1}-y_i=1$. Furthermore, we may assume that $-{3\over 2}\le x_1<{1\over 2}$ by a uni-modular transformation and assume that $G_i$ is primitive if it is a lattice vector by reduction.

If one of $G_2$, $G_3$ and $G_4$ is vertical, it can be easily deduced that
$${\rm vol}(P_8)\ge 6.\eqno (37)$$
For instance, when $G_3$ is vertical, we have $x_3-x_2\ge 1$, $x_4-x_5\ge 1$ and thus $x_3=x'_3=x_4\ge {3\over 2}$. Then, it can be deduced that
$${\rm vol}(P_8)\ge 4\cdot 2x_3-2(2x_3-1)=4x_3+2\ge 8,$$
which contradicts the assumption of (32).

Now, we assume that all $G_2$, $G_3$ and $G_4$ are not vertical.

\smallskip\noindent
{\bf Subcase 2.1.} {\it $y_3-y_2=2$ and ${\bf u}_2\not\in {1\over 2}\Lambda$}. Then ${\bf v}_3-{\bf v}_2=(k,2)$ is a lattice vector, where $k$ is a positive integer. On the other hand, it follows by the assumption $-{3\over 2}\le x_1<{1\over 2}$ that
$${\bf v}_5-{\bf v}_2=(x, 4),$$
where $-2<x\le 2$. Let $P$ denote the parallelogram with vertices ${\bf v}_1$, ${\bf v}_2$, ${\bf v}_5$ and ${\bf v}_6$, and let $T$ denote the triangle with vertices ${\bf v}_2$, ${\bf v}_3$ and ${\bf v}_5$, as shown by Figure 9.

\begin{figure}[!ht]
\centering
\includegraphics[scale=0.5]{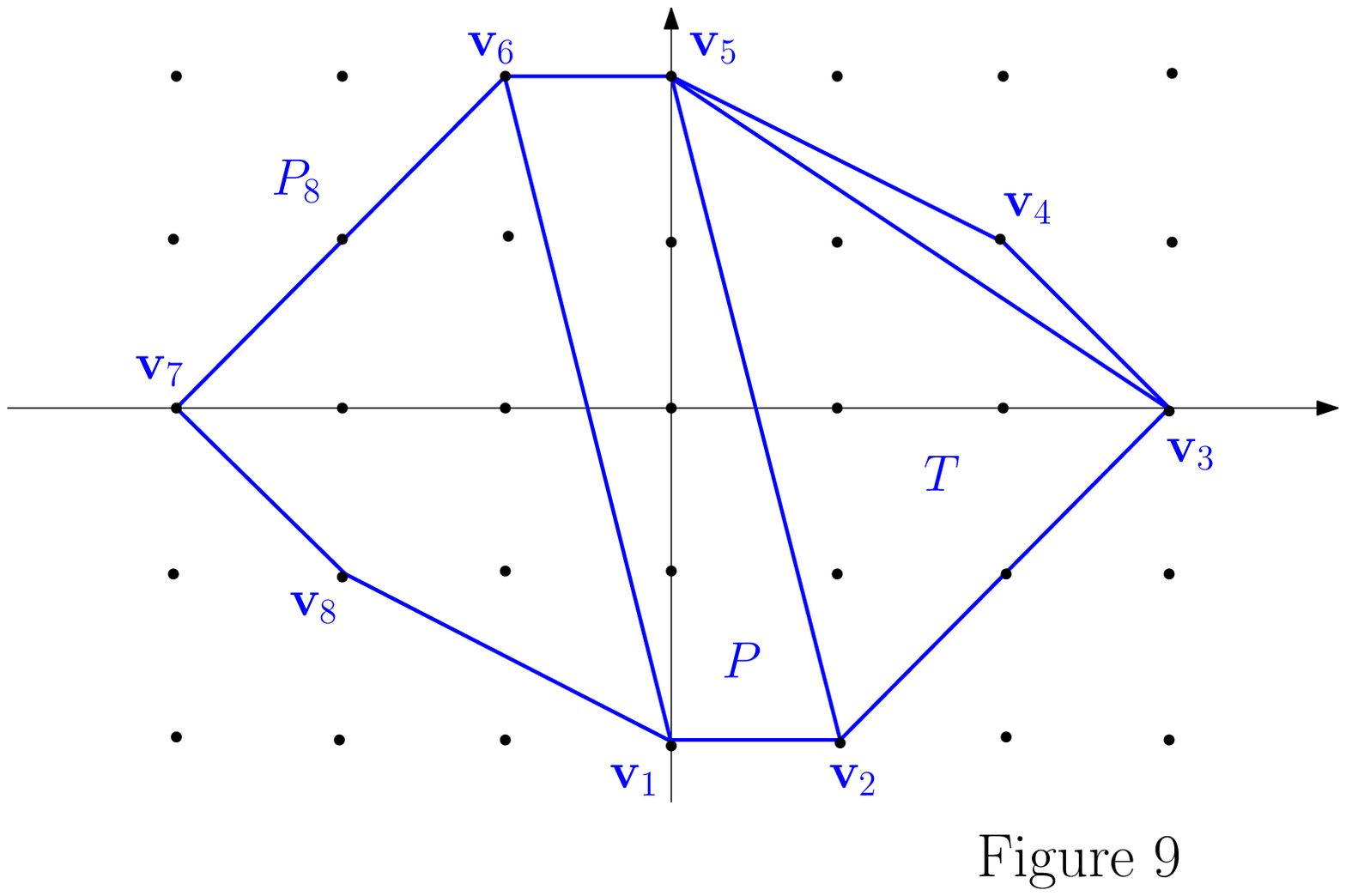}
\end{figure}

If $k\ge 2$, one can deduce
$${\rm vol}(T)={1\over 2}\left| \begin{array}{cc}
k&2\\
x&4
\end{array}\right|=2k-x\ge 2$$
and therefore
$${\rm vol}(P_8)> {\rm vol}(P)+2\cdot {\rm vol}(T)\ge 8,\eqno(38)$$
which contradicts the assumption of (32).

If $k=x_3-x_2=1$, $G_2\cap {1\over 2}\Lambda \not=\emptyset $ and ${\bf u}_2\not\in {1\over 2}\Lambda $, one can deduce that $x_2\in {1\over 4}\mathbb{Z}$ and therefore $x_1\in {1\over 4}\mathbb{Z}$. In fact, by checking all the eight cases $x_1=-{3\over 2},$ $-{5\over 4},$ $-1,$ $-{3\over 4},$ $-{1\over 2},$ $-{1\over 4},$ $0$ or ${1\over 4}$, it can be shown that there is no such octagon satisfying the conditions of Lemma 1.
For example, when $x_1={1\over 4}$, by convexity (as shown by Figure 10) the only candidate for ${\bf u}_3$ is ${\bf u}'_3=(2, {1\over 2})$ and the only candidates for ${\bf u}_4$ are ${\bf u}'_4=({1\over 2}, {3\over 2})$ and ${\bf u}^*_4=(1, {3\over 2})$. However, no octagon $P_8$ satisfying Lemma 1 can be constructed from these candidate midpoints.

\begin{figure}[!ht]
\centering
\includegraphics[scale=0.5]{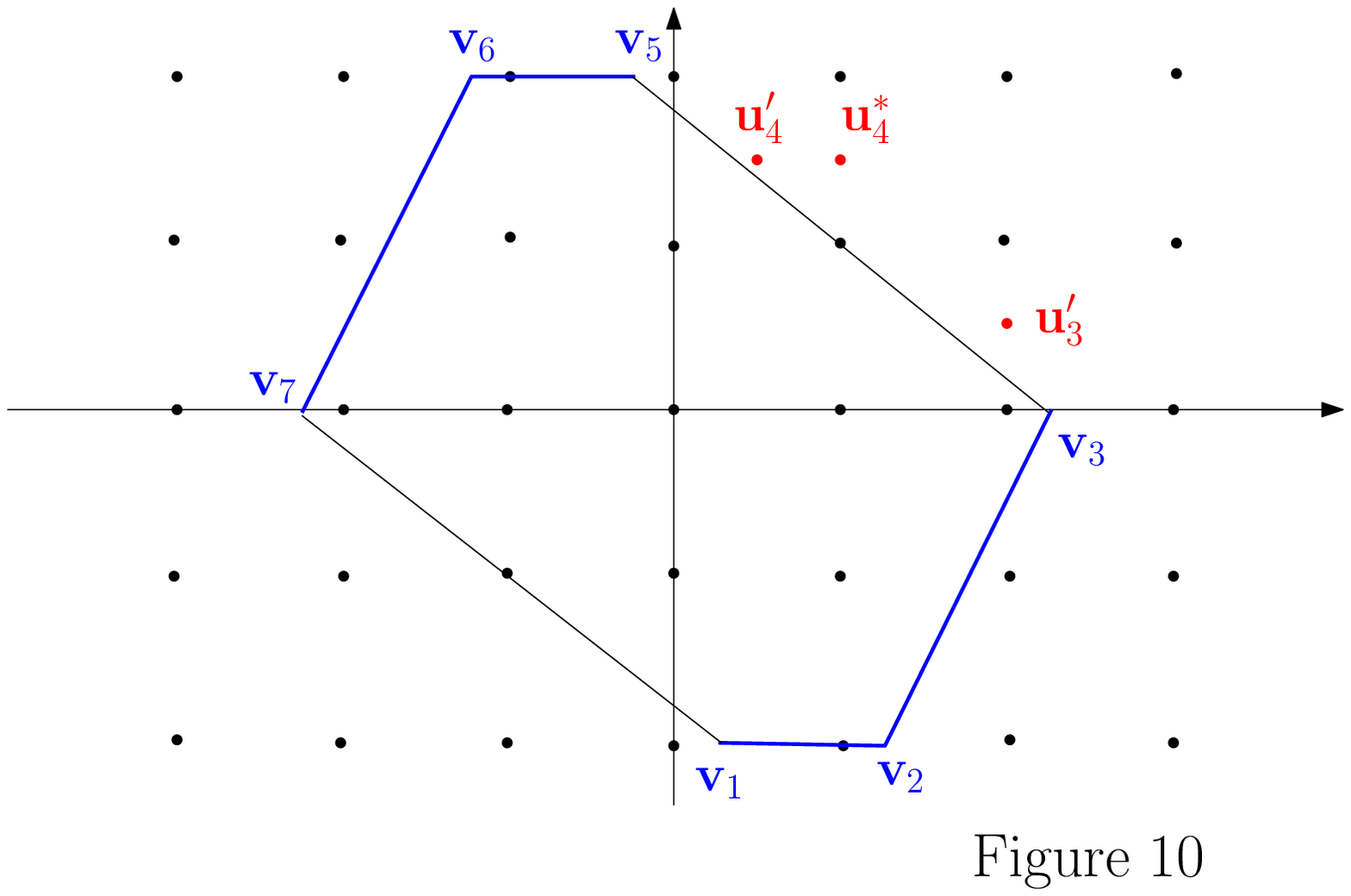}
\end{figure}

\noindent
{\bf Subcase 2.2.} {\it $y_4-y_3=2$ and ${\bf u}_3\not\in {1\over 2}\Lambda$.}  Then ${\bf v}_4-{\bf v}_3=(k,2)$ is a lattice vector, where $k$ is a positive integer (if it is negative, then make a reflection with respect to the $x$-axis). On the other hand, it follows by the assumption $-{3\over 2}\le x_1<{1\over 2}$ that
$${\bf v}_5-{\bf v}_2=(x, 4),$$
where $-2<x\le 2$. Let $P$ denote the parallelogram with vertices ${\bf v}_1$, ${\bf v}_2$, ${\bf v}_5$ and ${\bf v}_6$, and let $T$ denote the triangle with vertices ${\bf v}_2$, ${\bf v}'_3={\bf v}_2+({\bf v}_4-{\bf v}_3)$ and ${\bf v}_5$, as shown by Figure 11.

\begin{figure}[!ht]
\centering
\includegraphics[scale=0.5]{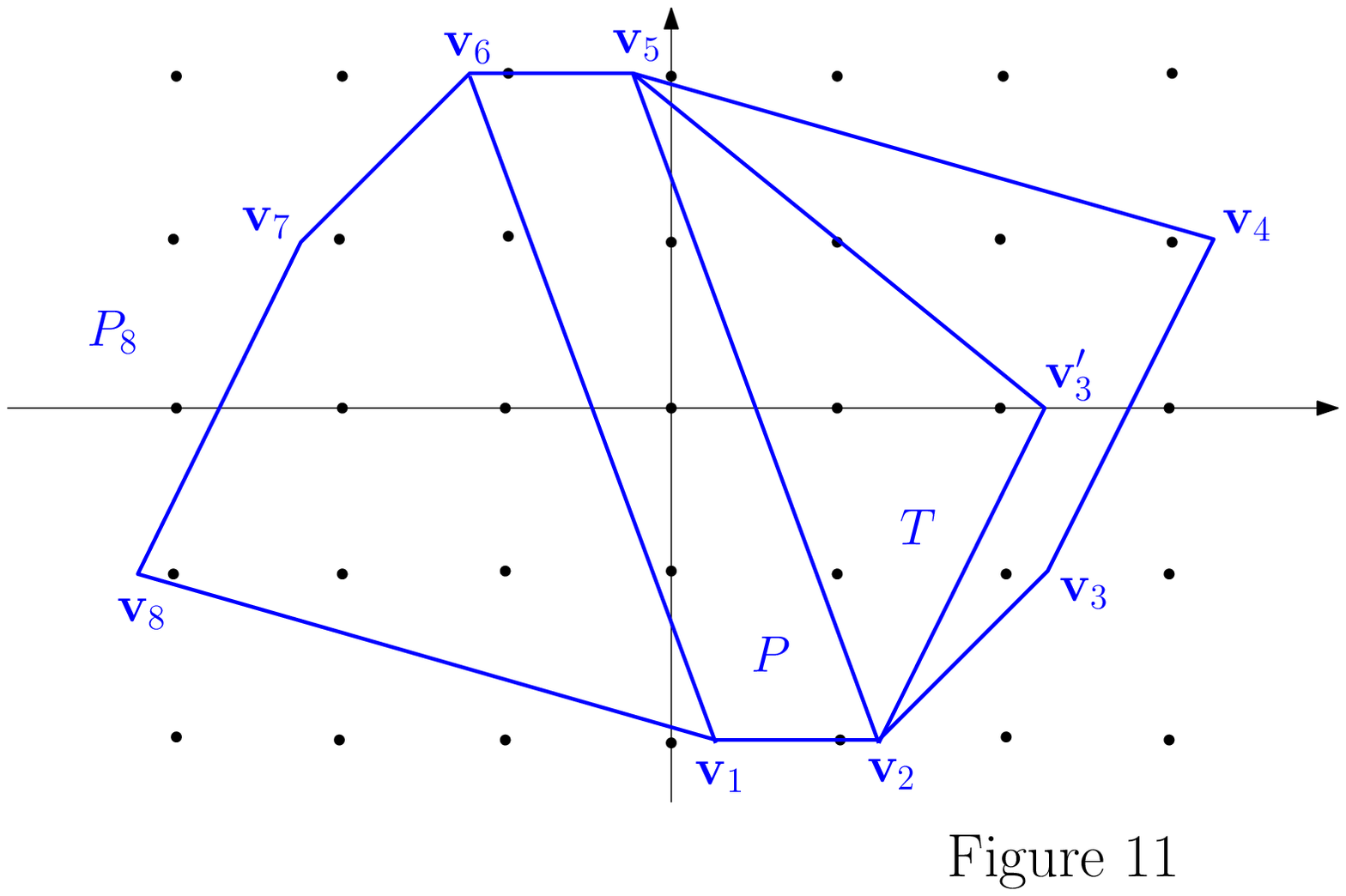}
\end{figure}

If $k\ge 2$, one can deduce
$${\rm vol}(T)={1\over 2}\left| \begin{array}{cc}
k&2\\
x&4
\end{array}\right|=2k-x\ge 2$$
and therefore
$${\rm vol}(P_8)> {\rm vol}(P)+2\cdot {\rm vol}(T)\ge 8,\eqno(39)$$
which contradicts the assumption of (32).

If $k=x_4-x_3=1$, $G_3\cap {1\over 2}\Lambda \not=\emptyset $ and ${\bf u}_3\not\in {1\over 2}\Lambda $, one can deduce that $x_3\in {1\over 4}\mathbb{Z}$ and therefore $x_1\in {1\over 4}\mathbb{Z}$. By checking all the eight cases $x_1=-{3\over 2},$ $-{5\over 4},$ $-1,$ $-{3\over 4},$ $-{1\over 2},$ $-{1\over 4},$ $0$ or ${1\over 4}$, it can be deduced that
$${\rm vol}(P_8)\ge 6.\eqno (40)$$
For example, when $x_1=-{3\over 2}$, we define ${\bf v}'_3=({3\over 2}, -1)$, ${\bf v}'_4=({5\over 2}, 1)$, ${\bf v}'_7=(-{3\over 2},1)$,
${\bf v}'_8=(-{5\over 2},-1)$, and define $P'_8$ to be the octagon with vertices ${\bf v}_1$, ${\bf v}_2$, ${\bf v}'_3$, ${\bf v}'_4$, ${\bf v}_5$, ${\bf v}_6$, ${\bf v}'_7$ and ${\bf v}'_8$, as shown by Figure 12. By shifting $G_3$ and $G_7$, one can deduce $P'_8\subseteq P_8$ and therefore
$${\rm vol}(P_8)\ge {\rm vol}(P'_8)=13.$$

\begin{figure}[!ht]
\centering
\includegraphics[scale=0.5]{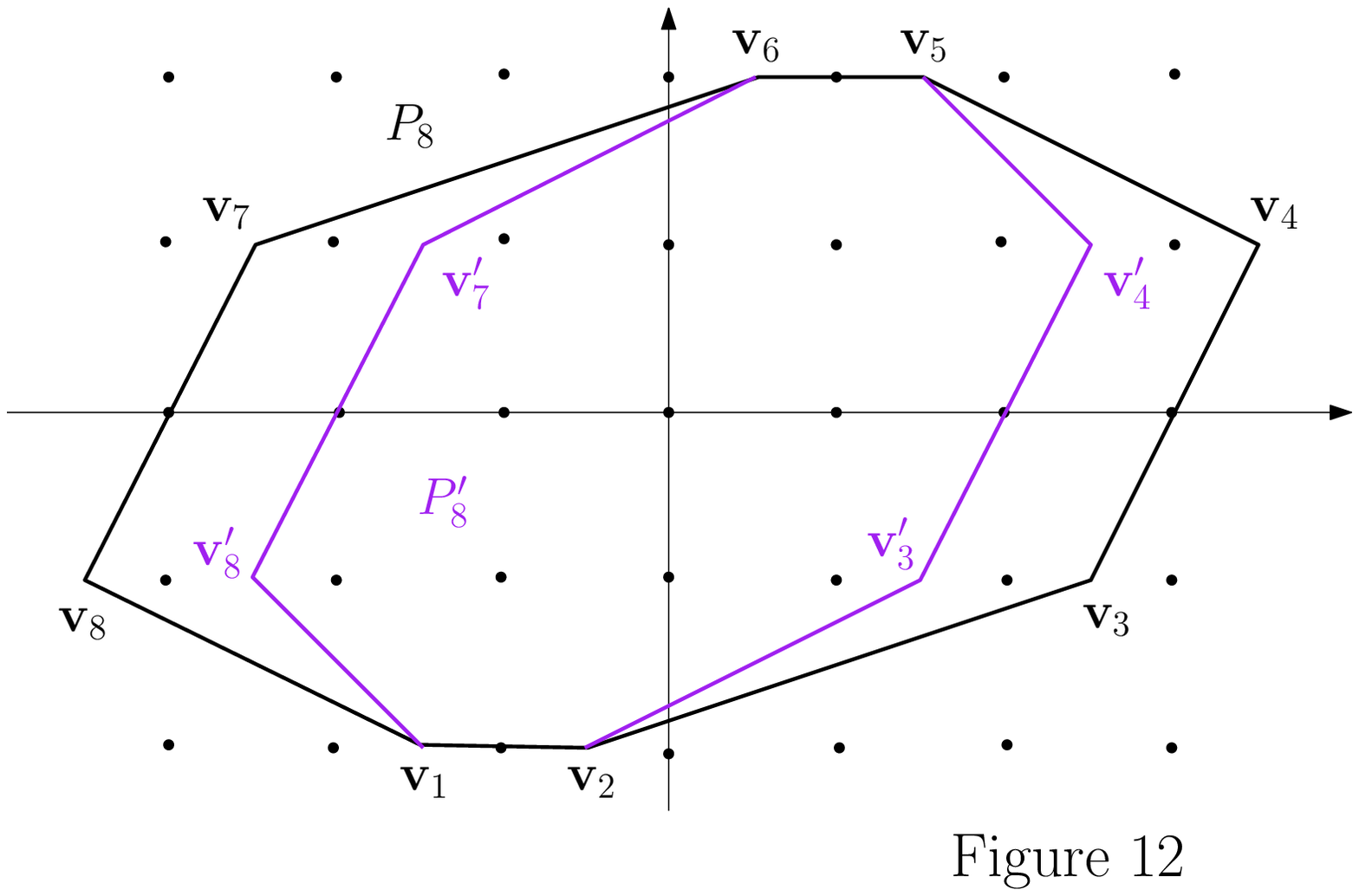}
\end{figure}

\noindent
{\bf Subcase 2.3.} {\it None of the three edges $G_2$, $G_3$ and $G_4$ is vertical and all ${\bf u}_2$, ${\bf u}_3$ and ${\bf u}_4$ belong to ${1\over 2}\Lambda $.} Then, it is sufficient to consider the following three situations.

\medskip\noindent
{\bf Subcase 2.3.1.}  $x'_3>\max\{ x'_2, x'_4\}$. Similar to Subcase 1.1, we get $x'_3\ge {3\over 2}$ and therefore
$${\rm vol}(P_8)\ge 4\cdot 2x'_3-2(2x'_3-1)=4x'_3+2\ge 8,\eqno (41)$$
which contradicts the assumption of (32).

\medskip
\noindent
{\bf Subcase 2.3.2.} $x'_2>\max \{ x'_3, x'_4\}$. If $x_3>x_2$, just like Subcase 1.2, one can get $x'_2\ge {3\over 2}$ and
$${\rm vol}(P_8)\ge 4\cdot 2x'_2-3(2x'_2-1)\ge 6,\eqno(42)$$
which contradicts the assumption of (32).

If $x_2>x_3$ and $y_3-y_2=1$, since $-{3\over 2}\le x_1<{1\over 2}$, ${\bf u}_2$ only can be $(1, -{3\over 2})$, $({1\over 2}, -{3\over 2})$,
$(0, -{3\over 2})$ or $(-{1\over 2}, -{3\over 2})$. Then it can be routinely checked that there is no convex octagon of this type satisfying Lemma 1.

If $x_2>x_3$ and $y_3-y_2=2$, since $-{3\over 2}\le x_1<{1\over 2}$, ${\bf u}_2$ only can be $(1, -1)$, $({1\over 2}, -1)$,
$(0, -1)$ or $(-{1\over 2}, -1)$. By checking these four cases, it can be shown that there is only one class of such convex octagons satisfying Lemma 1. Namely, the ones satisfying ${\bf u}_2=(1,-1)$, ${\bf u}_3=({1\over 2}, {1\over 2})$ and ${\bf u}_4=(0,{3\over 2})$, as shown in Figure 13. In other words, they are the octagons with vertices ${\bf v}_1=(\beta , -2)$, ${\bf v}_2=(1+\beta , -2)$, ${\bf v}_3=(1-\beta , 0)$, ${\bf v}_4=(\beta , 1)$, ${\bf v}_5=-{\bf v}_1$, ${\bf v}_6=-{\bf v}_2$, ${\bf v}_7=-{\bf v}_3$, ${\bf v}_8=-{\bf v}_4$, where ${1\over 4}<\beta <{1\over 3}$. Surprisingly, octagons of this type indeed satisfy
$${\rm vol}(P_8)=5.\eqno(43)$$
They are the octagons of second type listed in Lemma 9.

\begin{figure}[!ht]
\centering
\includegraphics[scale=0.5]{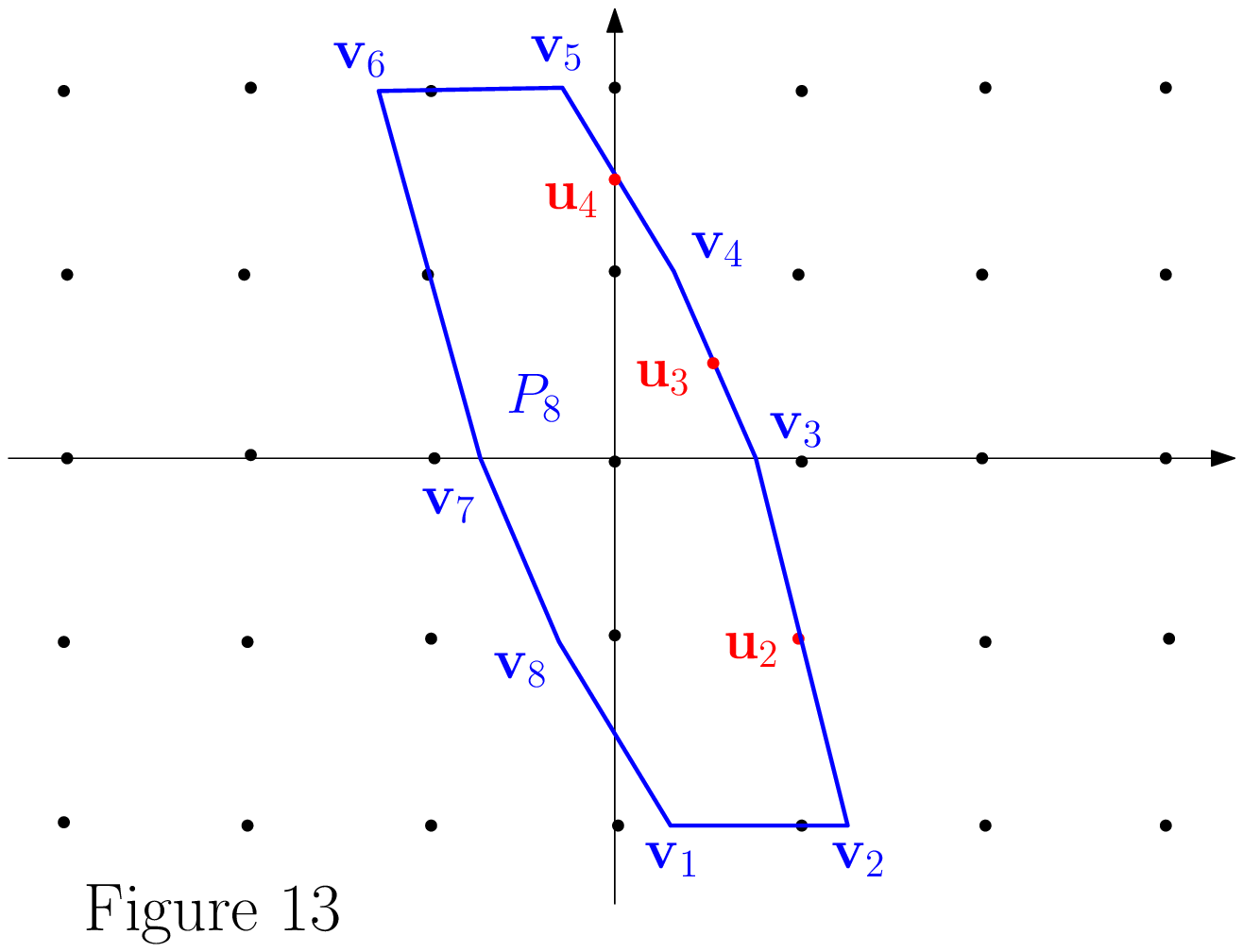}
\end{figure}

\medskip
\noindent
{\bf Subcase 2.3.3.} $x'_2=x'_3>x'_4$. Similar to Subcase 1.3, one can deduce $x'_2\ge 1$ and therefore
$${\rm vol}(P_8)\ge 4\cdot 2x'_3-2(2x'_3-1)=4x'_3+2\ge 6,\eqno (44)$$
which contradicts the assumption of (32).

As a conclusion of all these cases, Lemma 9 is proved. \hfill{$\Box$}

\vspace{0.6cm}
\noindent
{\Large\bf 4. Proofs of the Theorems}

\bigskip\noindent
{\bf Proof of Theorem 1.} Theorem 1 follows from Lemmas 5-9 immediately. \hfill{$\Box$}

\bigskip\noindent
{\bf Proof of Theorem 2.} Let $Q_{10}$ denote the convex decagon with vertices ${\bf u}_1=(0, 1)$, ${\bf u}_2=(1, 1)$, ${\bf u}_3=({3\over 2}, {1\over 2})$, ${\bf u}_4=({3\over 2}, 0)$, ${\bf u}_5=(1, -{1\over 2})$, ${\bf u}_6=-{\bf u}_1$, ${\bf u}_7=-{\bf u}_2$, ${\bf u}_8=-{\bf u}_3$, ${\bf u}_9=-{\bf u}_4$ and ${\bf u}_{10}=-{\bf u}_5$, let $L_i$ denote the straight line containing ${\bf u}_i$ and ${\bf u}_{i+1}$, where ${\bf u}_{10+i}={\bf u}_i$ and $L_{10+i}=L_i$, let ${\bf v}'_i$ denote the common point of $L_{i-2}$ and $L_i$, and let $T_i$ denote the triangle with vertices ${\bf v}'_i$, ${\bf u}_i$ and ${\bf u}_{i-1}$, as shown by Figure 14.

\begin{figure}[!ht]
\centering
\includegraphics[scale=0.5]{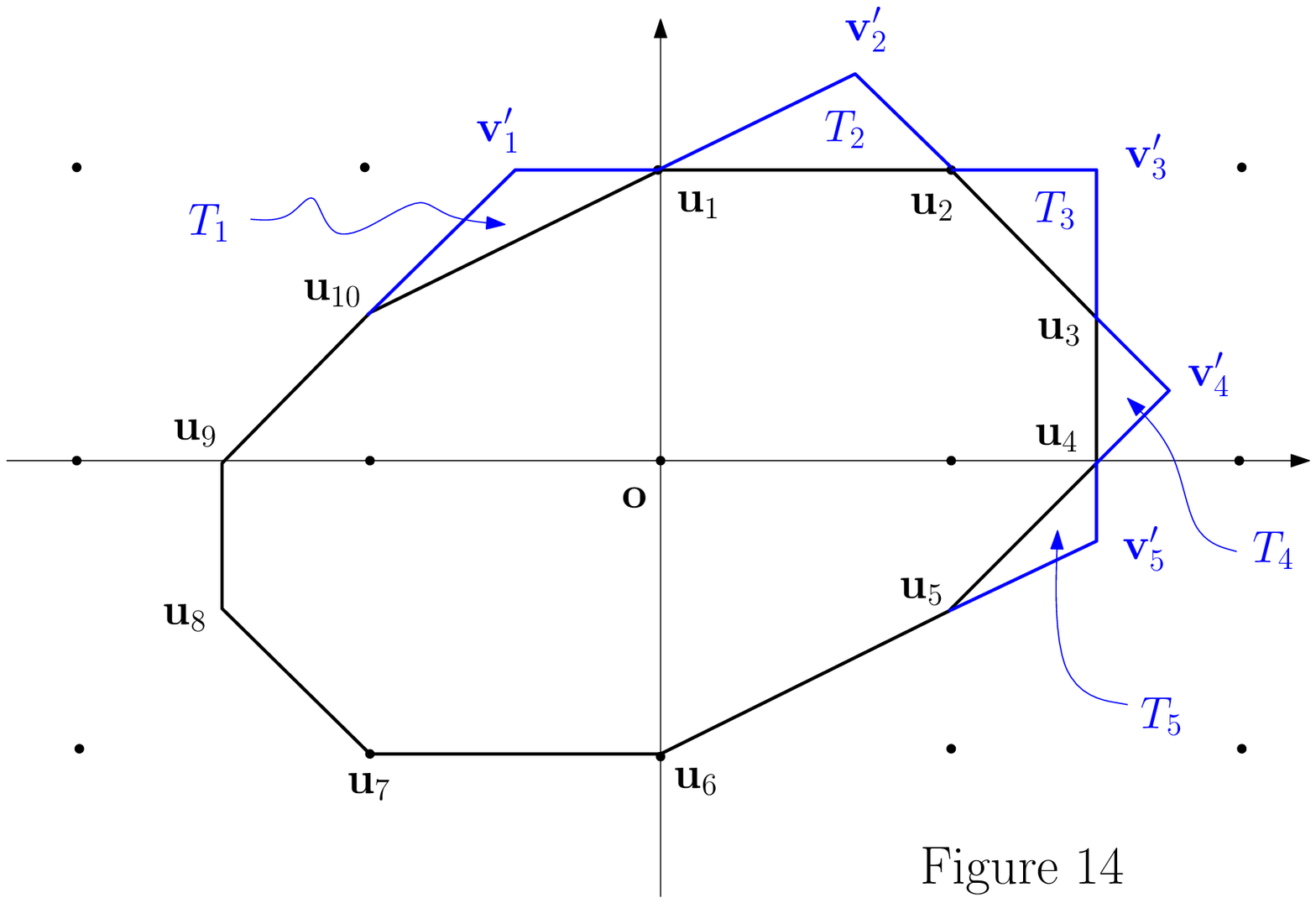}
\end{figure}

Assume that $P_{10}$ is a five-fold lattice tile with vertices ${\bf v}_1$, ${\bf v}_2$, $\ldots ,$ ${\bf v}_{10}$ satisfying
$${\bf v}_{i+1}-{\bf u}_i={\bf u}_i-{\bf v}_i$$
and therefore
$${\bf v}_{i+1}=2{\bf u}_i-{\bf v}_i, \eqno(45)$$
where ${\bf v}_{10+i}={\bf v}_i$. Apparently, it follows by convexity that
$${\bf v}_i\in {\rm int}(T_i), \quad i=1,\ 2,\ \ldots, \ 10.$$
In addition, by (45) we have
\begin{align*}
{\bf v}_5&\in {\rm int}(T_5),\\
{\bf v}_4&\in \bigl(2{\bf u}_4-{\rm int} (T_5)\bigr)\cap {\rm int}(T_4),\\
{\bf v}_3&\in \Bigl(2{\bf u}_3-\bigl(2{\bf u}_4-{\rm int}(T_5)\bigr)\cap {\rm int}(T_4)\Bigr)\cap {\rm int}(T_3)\\
& = \bigl(2({\bf u}_3-{\bf u}_4)+{\rm int}(T_5)\bigr)\cap \bigl(2{\bf u}_3-{\rm int}(T_4)\bigr)\cap {\rm int}(T_3),\\
{\bf v}_2&\in \Bigl(2{\bf u}_2-\bigl(2{\bf u}_3-\bigl(2{\bf u}_4-{\rm int}(T_5)\bigr)\cap {\rm int}(T_4)\bigr)\cap {\rm int}(T_3)\Bigr)\cap
{\rm int}(T_2)\\
& = \bigl(2({\bf u}_2-{\bf u}_3+{\bf u}_4)-{\rm int}(T_5)\bigr)\cap \bigl(2({\bf u}_2-{\bf u}_3)+{\rm int}(T_4)\bigr)\cap
 \bigl( 2{\bf u}_2-{\rm int}(T_3)\bigr)\cap {\rm int} (T_2),\\
{\bf v}_1&\in \Bigl( 2{\bf u}_1-\bigl(2{\bf u}_2-\bigl(2{\bf u}_3-\bigl(2{\bf u}_4-{\rm int}(T_5)\bigr)\cap {\rm int}(T_4)\bigr)\cap
{\rm int}(T_3)\bigr)\cap {\rm int}(T_2) \Bigr)\cap {\rm int}(T_1)\\
&= \bigl(2({\bf u}_1-{\bf u}_2+{\bf u}_3-{\bf u}_4)+{\rm int}(T_5)\bigr)\cap \bigl(2({\bf u}_1-{\bf u}_2+{\bf u}_3)-
{\rm int}(T_4)\bigr)\cap \\
& \quad\ \bigl( 2({\bf u}_1-{\bf u}_2)+{\rm int}(T_3)\bigr)\cap \bigl( 2{\bf u}_1-{\rm int} (T_2)\bigr) \cap {\rm int}(T_1).
\end{align*}
For convenience, we define
\begin{align*}
W&=\bigl(2({\bf u}_1-{\bf u}_2+{\bf u}_3-{\bf u}_4)+T_5\bigr)\cap \bigl(2({\bf u}_1-{\bf u}_2+{\bf u}_3)-T_4\bigr)\cap \\
&\quad\ \bigl( 2({\bf u}_1-{\bf u}_2)+T_3\bigr)\cap \bigl( 2{\bf u}_1-T_2\bigr) \cap T_1.
\end{align*}

On the other hand, whenever we take
$${\bf v}_1\in {\rm int}(W)$$
and define ${\bf v}_i$ successively by (45), the inverse of the above process and Lemma 3 guarantee that
$${\bf v}_i\in {\rm int}(T_i)$$
holds for all $i=1,$ $2,$ $\ldots,$ $10$. Therefore, by Lemma 1 the decagon with them as its vertices is indeed a five-fold lattice tile.

By routine and detailed computation, it can be deduced from its definition that $W$ is a quadrilateral with vertices ${\bf w}_1=(-{1\over 2},1)$, ${\bf w}_2=(-{1\over 2},{3\over 4})$, ${\bf w}_3=(-{2\over 3}, {2\over 3})$ and ${\bf w}_4=(-{3\over 4},{3\over 4})$. Theorem 2 is proved. \hfill{$\Box$}

\vspace{0.6cm}\noindent
{\bf Acknowledgements.} This work is supported by 973 Program 2013CB834201.

\bibliographystyle{amsplain}

\vspace{0.6cm}
\noindent
Chuanming Zong, Center for Applied Mathematics, Tianjin University, Tianjin 300072, China

\noindent
Email: cmzong@math.pku.edu.cn

\end{document}